\documentclass[a4paper]{amsart}

\usepackage{amsmath,amssymb,amsfonts}
\usepackage{exscale,url,doi,graphicx,psfrag,xcolor}
\usepackage[ruled,vlined,algosection]{algorithm2e}

\definecolor{color1}{rgb}{0,0,1}
\hypersetup{pdfborder={1 0 0}}

\def\diag{\ensuremath{\mathrm{diag}}}
\def\la{\ensuremath{\lambda}}
\def\RR{\ensuremath{\mathbb R}}
\def\CC{\ensuremath{\mathbb C}}
\def\spanl{\ensuremath{\mathrm{span}}}

\newtheorem{theorem}{Theorem}[section]
\newtheorem{lemma}[theorem]{Lemma}

\theoremstyle{definition}

\newtheorem{remark}[theorem]{Remark}
\numberwithin{equation}{section}

\begin{document}

\title[Individual Ritz values in BPG]
{Convergence rates of individual Ritz values\\
 in block preconditioned gradient-type eigensolvers}

\author[M.~Zhou]{Ming Zhou}
\address{Universit\"at Rostock, Institut f\"ur Mathematik, 
        Ulmenstra{\ss}e 69, 18055 Rostock, Germany}
\email{ming.zhou at uni-rostock (dot) de}

\author[K.~Neymeyr]{Klaus Neymeyr}
\address{Universit\"at Rostock, Institut f\"ur Mathematik, 
        Ulmenstra{\ss}e 69, 18055 Rostock, Germany}
\email{klaus.neymeyr at uni-rostock (dot) de}


\keywords{preconditioned subspace eigensolvers, Ritz values, cluster robustness.
 \hfill June 1, 2022}

\begin{abstract}
Many popular eigensolvers for large and sparse Hermitian matrices
or matrix pairs can be interpreted as accelerated
block preconditioned gradient (BPG) iterations in order to
analyze their convergence behavior by composing known estimates.
An important feature of BPG is the cluster robustness, i.e.,
reasonable performance for computing clustered eigenvalues
is ensured by a sufficiently large block size.
This feature can easily be explained
for exact-inverse (exact shift-inverse) preconditioning
by adapting classical estimates on nonpreconditioned eigensolvers,
whereas the existing results for more general preconditioning
are still improvable. We expect to extend certain sharp estimates
for the corresponding vector iterations to BPG where proper bounds
of convergence rates of individual Ritz values are to be derived.
Such an extension has been achieved for BPG with fixed step sizes in
[Math. Comp. 88 (2019), 2737--2765]. The present paper
deals with the more practical case that the step sizes
are implicitly optimized by the Rayleigh-Ritz method.
Our new estimates improve some previous ones
in view of concise and more flexible bounds.
\end{abstract}

\maketitle
\pagestyle{myheadings}
\thispagestyle{plain}

\section{Introduction}

Solving eigenvalue problems of large and sparse Hermitian matrices
or matrix pairs are of practical importance in various applications.
Proper iterative methods with vectors or subspaces allow determining
the desired eigenpairs with reasonable effort \cite{TEV2000,KNN2003e}.
The convergence behavior of such eigensolvers depends on
the distribution of relevant eigenvalues as well as
certain methodical characteristics including preconditioners
and block sizes (dimensions of iterates).

As a simple example, we first consider the computation of the smallest eigenvalues
of a symmetric positive definite matrix $A\in\RR^{n \times n}$
by the preconditioned subspace iteration
\begin{equation}\label{bpinvit}
 X^{(\ell+1)} \quad\xleftarrow{\mathrm{RR}[s]}\quad \spanl\{X^{(\ell)}-TR^{(\ell)}\}.
\end{equation}
Therein $s$ denotes the block size.
The current iterate $X^{(\ell)}\in\RR^{n \times s}$ is assumed to have full rank
and consists of orthonormal Ritz vectors of $A$ in the subspace $\spanl\{X^{(\ell)}\}$.
The corresponding residuals form the block residual
$R^{(\ell)}=AX^{(\ell)}-X^{(\ell)}\Theta^{(\ell)}$
with the diagonal matrix $\Theta^{(\ell)}\in\RR^{s \times s}$
containing Ritz values. The term $TR^{(\ell)}$ can be determined
by using an incomplete factorization of $A$ or
approximately solving a block linear system of the form $AE=R^{(\ell)}$.
The underlying matrix $T$ is called preconditioner
and represents an approximate inverse of $A$ for which the condition
\begin{equation}\label{prec}
 \|I-TA\|_A\le\gamma<1
\end{equation}
with the $n{\times}n$ identity matrix $I$ ensures that
the trial subspace $\mathcal{U}^{(\ell)}=\spanl\{X^{(\ell)}-TR^{(\ell)}\}$
has dimension $s$ according to \cite[Lemma 3.1]{NEY2000ps}.
Finally, the Rayleigh-Ritz procedure $\mathrm{RR}[s]$
extracts orthonormal Ritz vectors from $\mathcal{U}^{(\ell)}$
and builds with them the next iterate $X^{(\ell+1)}$.
This elementary eigensolver is actually a block preconditioned gradient
(BPG) iteration since the columns of $R^{(\ell)}$
are collinear with the gradient vectors of the Rayleigh quotient
\[\la:\RR^n{\setminus}\{0\}\to\RR, \quad
 \la(x)=\frac{x^TAx}{x^Tx}\]
associated with the columns of $X^{(\ell)}$.
For computing the first $t$ eigenvalues of $A$
concerning the eigenvalue arrangement $\la_1\le\cdots\le\la_n$,
one can basically set a sufficiently large block size $s$ for ensuring
$\la_t\ll\la_{s+1}$ according to the well-known convergence theory
\cite{PAR1980} for the subspace iteration
$X^{(\ell+1)} \ \xleftarrow{\mathrm{RR}[s]} \ \spanl\{A^{-1}X^{(\ell)}\}$
(i.e. the block power method for $A^{-1}$)
which coincides with the special form of \eqref{bpinvit} for $T=A^{-1}$.
Extending the trial subspace of \eqref{bpinvit}
leads to more efficient eigensolvers such as
\begin{equation}\label{bpg}
 X^{(\ell+1)} \quad\xleftarrow{\mathrm{RR}[s]}\quad \spanl\{X^{(\ell)},TR^{(\ell)}\}
\end{equation}
which can be interpreted as a BPG iteration with optimized step sizes, and
\begin{equation}\label{lobpcg}
 X^{(\ell+1)} \quad\xleftarrow{\mathrm{RR}[s]}\quad \spanl\{X^{(\ell-1)},X^{(\ell)},TR^{(\ell)}\}
\end{equation}
which corresponds to the locally optimal block
preconditioned conjugate gradient (LOBPCG) method \cite{KNY2001}.
Moreover, it is advantageous to implement these eigensolvers in combination with deflation
due to the different convergence rates of individual Ritz values.
Deriving sharp bounds of these convergence rates
is challenging for advanced eigensolvers. We review here some known results
for the iterations \eqref{bpinvit} and \eqref{bpg}.

\subsection{Known results}

The convergence behavior of \eqref{bpinvit}
can be analyzed as in \cite[Section 2]{BPK1996} in terms of
the eigenvalue arrangement $\la_1\le\cdots\le\la_n$ of $A$,
the quality parameter $\gamma$ from \eqref{prec} for the preconditioner $T$,
and the block size $s$. The resulting estimate
for the $i$th Ritz value for an index $i\in\{1,\ldots,s\}$
provides a bound which essentially depends on the convergence factor
\begin{equation}\label{convf}
 \gamma+(1-\gamma)\la_i/\la_{s+1}.
\end{equation}
Its special form $\la_i/\la_{s+1}$ for $\gamma=0$, i.e., for $T=A^{-1}$,
also appears in a classical angle estimate
for the block power method for $A^{-1}$ in \cite{PAR1980}.
However, the analysis in \cite{BPK1996} requires
a technical assumption on certain angles between
the initial subspace $\spanl\{X^{(0)}\}$ and the eigenvectors associated with
the eigenvalues $\la_1,\ldots,\la_s$. The gap $\la_{i+1}-\la_i$ for each
$i\in\{1,\ldots,s\}$ has to be sufficiently large for making the assumption
practically reasonable. Although the convergence factor
\eqref{convf} is suitable for indicating the cluster robustness of \eqref{bpinvit},
i.e., fast convergence for $i \ll s$ despite $\la_{i+1}-\la_i\approx0$,
the assumption limits the applicability.

A more flexible and concise estimate for \eqref{bpinvit}
can be derived by \cite[Section 5]{KNN2003}.
Therein an eigenvalue interval $(\la_j,\la_{j+1})$ with $j \ge i$ is used for locating
the $i$th Ritz value in the current subspace iterate, and
the distance ratio $(*-\la_j)/(\la_{j+1}-*)$ serves as a convergence measure.
The corresponding convergence factor reads $\gamma+(1-\gamma)\la_j/\la_{j+1}$.
In particular, if the $i$th Ritz value reaches the interval $(\la_i,\la_{i+1})$,
one gets the special form $\gamma+(1-\gamma)\la_i/\la_{i+1}$ which is less accurate
in comparison to \eqref{convf} but still reasonable for sufficiently large $\la_{i+1}-\la_i$.
Moreover, this result cannot be refined as \eqref{convf} without further modifications
since its theoretical sharpness can be verified by certain special iterates.

Our recent result in \cite{NEZ2019} uses a larger interval
for the Ritz value location, namely, $(\la_{j-s+i},\la_{j+1})$ with $j \ge s$.
By defining an alternative quality parameter $\widetilde{\gamma}$
for the preconditioner $T$ concerning a geometric interpretation
based on \cite{NEY2000ps}, we have achieved the convergence factor
$\widetilde{\gamma}+(1-\widetilde{\gamma})\la_{j-s+i}/\la_{j+1}$
with respect to $(*-\la_{j-s+i})/(\la_{j+1}-*)$.
The final phase of \eqref{bpinvit} is characterized by $j=s$ where the distance ratio
is simply $(*-\la_i)/(\la_{s+1}-*)$ and the convergence factor is specialized as
$\widetilde{\gamma}+(1-\widetilde{\gamma})\la_i/\la_{s+1}$.
This is comparable with \eqref{convf}, and can
reasonably describe the cluster robustness
since the technical assumption used in \cite{BPK1996} is avoided.

The above estimates for \eqref{bpinvit} also provide preliminary bounds for
the accelerated iterations \eqref{bpg} and \eqref{lobpcg}. More direct bounds
for \eqref{bpg} have been presented in \cite{OVT2006} in terms of
sums of Ritz value errors by generalizing some arguments
from \cite{SAM1958,OVT2006s} concerning vectorial gradient iterations,
and in \cite{NEZ2014} by upgrading the analysis
from \cite{KNN2003} and adapting a sharp estimate from \cite{NEY2012}
for the single-vector version of \eqref{bpg}. In comparison to
the desirable convergence factor
\begin{equation}\label{convf1}
 \frac{\kappa+\gamma(2-\kappa)}{(2-\kappa)+\gamma\kappa}
 \quad\mbox{with}\quad
 \kappa=\frac{\la_i(\la_n-\la_{s+1})}{\la_{s+1}(\la_n-\la_i)}
\end{equation}
which can easily be shown to improve \eqref{convf},
the result from \cite{OVT2006} contains some additional technical terms
and only asymptotically indicates \eqref{convf1}, whereas
the result from \cite{NEZ2014} gives a suboptimal form of \eqref{convf1}
with the index $i{\,+\,}1$ instead of $s{\,+\,}1$.
These results for \eqref{bpg} need to be improved
in order to build a proper basement for the convergence analysis
of more complicated eigensolvers including LOBPCG.

\subsection{Aim and overview}

Our goal is to derive concise Ritz value estimates
containing convergence factors like \eqref{convf1}
for interpreting the cluster robustness of the BPG iteration \eqref{bpg}.
As the ratio $\la_i/\la_{s+1}$ is a decisive term in \eqref{convf1},
a fundamental idea is to skip the eigenvalues
$\la_{i+1},\ldots,\la_s$ (or $\la_{j-s+i+1},\ldots,\la_j$
in a more general case with $j \ge s$) by utilizing
certain auxiliary subspaces which are orthogonal (and $A$-orthogonal) to
the associated eigenvectors $x_{i+1},\ldots,x_s$.

This idea arises from the analysis of an abstract block iteration
by Knyazev \cite{KNY1987}, and has been adapted
to the preconditioned subspace iteration \eqref{bpinvit} in \cite{NEZ2019}.
By observing a partial iteration of \eqref{bpinvit}
within the orthogonal complement of $\spanl\{x_{i+1},\ldots,x_s\}$,
some Ritz vectors in two successive subspace iterates are compared
in a geometric way similarly to \cite{NEY2000ps}
for constructing a perturbed inverse vector iteration.
The corresponding perturbation parameter $\widetilde{\gamma}$
can be used as an alternative quality parameter of preconditioning
in the further analysis. We note that this approach depends on the fact
that the next subspace iterate in \eqref{bpinvit}
is simply the current trial subspace, i.e.,
$\spanl\{X^{(\ell+1)}\}=\spanl\{X^{(\ell)}-TR^{(\ell)}\}$.
Thus a direct comparison between Ritz vectors is enabled.

In contrast, the BPG iteration \eqref{bpg} cannot be described by
an equality formula since the next subspace iterate is only a subset
of the current trial subspace. The more complicated relation
between Ritz vectors therein is analyzed in \cite{NEZ2014}
using Sion's-minimax theorem via certain basis matrices
instead of Ritz vectors. Consequently,  for analyzing
the cluster robustness of \eqref{bpg}, we tend to update
the above construction of an alternative quality parameter
of preconditioning by means of basis matrices or subspaces.

For the sake of generality, we follow the introduction of
the LOBPCG method in \cite{KNY2001} and reformulate
\eqref{bpg} for the generalized eigenvalue problem
\begin{equation}\label{evp}
 Mv=\mu Av, \quad M,A\in\CC^{n \times n}\ \mbox{Hermitian}, \quad
 A\ \mbox{positive definite}
\end{equation}
where the target eigenvalues of the matrix pair $(M,A)$ are the largest ones.
Some conversions between \eqref{evp} and practical eigenvalue problems
are introduced in Section 2 together with a simple representation
of the investigated iteration which does not limit the generality.
Section 3 provides some auxiliary terms and intermediate arguments
based on the analysis of an abstract block iteration from \cite{KNY1987}
and the analysis of the preconditioned subspace iteration \eqref{bpinvit}
from \cite{NEZ2019}. A subspace-oriented interpretation
of preconditioning in partial iterations of the BPG iteration \eqref{bpg}
leads to multi-step estimates in Section 4 reflecting the cluster robustness.
We additionally discuss the possibility of deriving estimates
under classical conditions like \eqref{prec}. Numerical experiments
for illustrating the new results are given in Section 5.

\section{Preliminaries}

The generalized eigenvalue problem \eqref{evp} can be used as a common form
of several practical eigenvalue problems, e.g., computing a subset of the spectrum
of a self-adjoint elliptic partial differential operator together with
the associated eigenfunctions. Therein proper discretizations produce
the standard eigenvalue problem $Lu=\la u$
of an Hermitian matrix $L\in\CC^{n \times n}$
or the generalized eigenvalue problem
\begin{equation}\label{evp1}
 Lu=\la Su, \quad L,S\in\CC^{n \times n}\ \mbox{Hermitian}, \quad
 S\ \mbox{positive definite}
\end{equation}
which formally includes $Lu=\la u$
by setting $S$ as the $n{\times}n$ identity matrix $I$.

If the target eigenvalues of the matrix pair $(L,S)$ are the smallest ones,
we can transform \eqref{evp1} as $(L-\sigma S)u=(\la-\sigma)Su$
with a sufficiently small shift $\sigma$ such that the matrix
$\widetilde{L}=L-\sigma S$ is positive definite. The shifted problem
corresponds to \eqref{evp} for $M=S$, $A=\widetilde{L}$ and $\mu=(\la-\sigma)^{-1}$.

If some interior eigenvalues of $(L,S)$ are first to be determined, a similar shifted problem
with an indefinite and invertible $\widetilde{L}$ can be established. This cannot directly
be covered by \eqref{evp}. Instead, we can consider the equivalent problem
$(\widetilde{L}S^{-1}\widetilde{L})u=(\la-\sigma)\widetilde{L}u$ as a special form of \eqref{evp} with
$M=\pm\widetilde{L}$, $A=\widetilde{L}S^{-1}\widetilde{L}$ and $\mu=\pm(\la-\sigma)^{-1}$.

Some further specializations of \eqref{evp} refer to Hermitian definite matrix pencils
\cite{KPS2014} and the linear response eigenvalue problem \cite{BEL2017}.

\subsection{Considered iteration}

We modify the iteration \eqref{bpg} in order to compute
the largest eigenvalues of $(M,A)$ from \eqref{evp}. With the block size $s$,
the current iterate $V^{(\ell)}\in\CC^{n \times s}$ has full rank
and consists of $A$-orthonormal Ritz vectors of $(M,A)$
in the subspace $\spanl\{V^{(\ell)}\}$. The associated block residual
reads $R_V^{(\ell)}=MV^{(\ell)}-AV^{(\ell)}\Theta_V^{(\ell)}$
with the diagonal Ritz value matrix
$\Theta_V^{(\ell)}={V^{(\ell)}}^*MV^{(\ell)}\in\RR^{s \times s}$.
An approximate solution of the block linear system $AE=R_V^{(\ell)}$
is denoted by $\widetilde{T}R_V^{(\ell)}$
with a Hermitian positive definite preconditioner $\widetilde{T}$
which is an approximate inverse of $A$. By using the smallest eigenvalue $\alpha$
and the largest eigenvalue $\beta$ of the matrix product $\widetilde{T}A$
(or $A^{1/2}\widetilde{T}A^{1/2}$) which are both positive, it holds that
\begin{equation}\label{prec1}
 \|I-\omega \widetilde{T}A\|_A\le\gamma
 \quad\mbox{with}\quad
 \omega=\frac{2}{\beta+\alpha}
 \quad\mbox{and}\quad
 \gamma=\frac{\beta-\alpha}{\beta+\alpha}<1.
\end{equation}
This condition is a more natural form of \eqref{prec} concerning an arbitrary
Hermitian positive definite $\widetilde{T}$ and additional scaling.
The trial subspace $\mathcal{U}_V^{(\ell)}=\spanl\{V^{(\ell)},\widetilde{T}R_V^{(\ell)}\}$
evidently has at least dimension $s$.
The next iterate $V^{(\ell+1)}$ is constructed by
$A$-orthonormal Ritz vectors of $(M,A)$ in $\mathcal{U}_V^{(\ell)}$
associated with the $s$ largest Ritz values.
We denote by $\mathrm{RR}[M,A,s]$ the underlying Rayleigh-Ritz procedure.
Then the modified version of \eqref{bpg} is represented by
\begin{equation}\label{bpg1}
 V^{(\ell+1)} \quad\xleftarrow{\mathrm{RR}[M,A,s]}\quad
 \spanl\{V^{(\ell)},\widetilde{T}R_V^{(\ell)}\}.
\end{equation}

The special form of \eqref{bpg1} for $M=I$ is equivalent to \eqref{bpg}.
Therein all Ritz values are positive
so that $\Theta_V^{(\ell)}$ is positive definite. By using its square root matrix
$C=(\Theta_V^{(\ell)})^{1/2}$, one can construct the iterate
$X^{(\ell)}=V^{(\ell)}C^{-1}$ for \eqref{bpg} due to the properties
\[{X^{(\ell)}}^*X^{(\ell)}
 ={X^{(\ell)}}^*MX^{(\ell)}
 =C^{-1}{V^{(\ell)}}^*MV^{(\ell)}C^{-1}
 =C^{-1}\Theta_V^{(\ell)}C^{-1}=I_s,\]
\[{X^{(\ell)}}^*AX^{(\ell)}
 =C^{-1}{V^{(\ell)}}^*AV^{(\ell)}C^{-1}
 =C^{-2}=(\Theta_V^{(\ell)})^{-1},\]
i.e., the columns of $X^{(\ell)}$ are orthonormal Ritz vectors of $A$ in $\spanl\{X^{(\ell)}\}$,
and the corresponding Ritz values are contained in $\Theta^{(\ell)}=(\Theta_V^{(\ell)})^{-1}$.
In addition, the relation
\[\begin{split}
 R^{(\ell)}&=AX^{(\ell)}-X^{(\ell)}\Theta^{(\ell)}
 =AV^{(\ell)}C^{-1}-V^{(\ell)}C^{-1}(\Theta_V^{(\ell)})^{-1}\\[1ex]
 &=-(MV^{(\ell)}-AV^{(\ell)}\Theta_V^{(\ell)})C^{-3}
 =-R_V^{(\ell)}C^{-3}
\end{split}\]
leads to the subspace equality
$\spanl\{X^{(\ell)},TR^{(\ell)}\}=\spanl\{V^{(\ell)},\widetilde{T}R_V^{(\ell)}\}$
for $T=\omega\widetilde{T}$ concerning the conditions \eqref{prec} and \eqref{prec1}.

Furthermore, if \eqref{bpg1} is applied to the practical problem \eqref{evp1},
determining $\widetilde{T}R_V^{(\ell)}$ refers to solving the block linear system
$AE=R_V^{(\ell)}$ for $A=\widetilde{L}$ or $A=\widetilde{L}S^{-1}\widetilde{L}$.
The latter case can be implemented by solving two systems for $\widetilde{L}$ successively.

\subsection{Convergence measure and simplified notation}

We denote by $\mu_i$ the $i$th largest eigenvalue of $(M,A)$ from \eqref{evp},
i.e., the eigenvalues are arranged as $\mu_1\ge\cdots\ge\mu_n$.
With the Ritz values $\theta_1^{(\ell)}\ge\cdots\ge\theta_s^{(\ell)}$
of $(M,A)$ in the subspace $\spanl\{V^{(\ell)}\}$, we set
$\Theta_V^{(\ell)}=\diag(\theta_1^{(\ell)},\ldots,\theta_s^{(\ell)})$.
We measure the convergence of $\theta_i^{(\ell)}$ by the distance ratio
$(\mu_{j-s+i}-*)/(*-\mu_{j+1})$ with a certain index $j \ge s$.

For an arbitrary shift $\sigma<\mu_n$, the iteration \eqref{bpg1} and
the above convergence measure are invariant for the substitution
\begin{equation}\label{subst}
 M \ \leftrightarrow \ M-\sigma A, \qquad
 \mu_i \ \leftrightarrow \ \mu_i-\sigma, \qquad
 \theta_i^{(\ell)} \ \leftrightarrow \ \theta_i^{(\ell)}-\sigma.
\end{equation}
Therefore assuming that $M$ is positive definite does not limit the generality.

In addition, since the condition \eqref{prec1} is formulated with respect to
the inner product induced by $A$, we can simplify the notation
of matrices and vectors by using the representations
\begin{equation}\label{rpr}
 H=A^{-1/2}MA^{-1/2}, \qquad y=A^{1/2}v, \qquad Y=A^{1/2}V,
 \qquad N=A^{1/2}(\omega \widetilde{T})A^{1/2}
\end{equation}
as in \cite[Subsection 1.2]{NEZ2019}.
Then \eqref{prec1} turns into
\begin{equation}\label{prec2}
 \|I-N\|_2\le\gamma<1,
\end{equation} 
and \eqref{bpg1} is equivalent to
\begin{equation}\label{bpg2}
 Y^{(\ell+1)} \quad\xleftarrow{\mathrm{RR}[H,s]}\quad
 \spanl\{Y^{(\ell)},N(HY^{(\ell)}-Y^{(\ell)}\Theta_Y^{(\ell)})\}
\end{equation}
where $\Theta_Y^{(\ell)}=\Theta_V^{(\ell)}$.
The notation of eigenvalues and Ritz values remains unchanged.

\begin{remark}\label{cas}
For analyzing the convergence behavior of the BPG iteration \eqref{bpg1},
we only need to observe the accompanying iteration \eqref{bpg2}
for two Hermitian positive definite matrices:
$H$ with the arranged eigenvalues $\mu_1\ge\cdots\ge\mu_n$,
and $N$ satisfying \eqref{prec2}.
Therein the current iterate $Y^{(\ell)}\in\CC^{n \times s}$ has full rank
and its columns are orthonormal Ritz vectors of $H$ in $\spanl\{Y^{(\ell)}\}$
associated with the arranged Ritz values
$\theta_1^{(\ell)}\ge\cdots\ge\theta_s^{(\ell)}$,
also contained in the diagonal matrix $\Theta_Y^{(\ell)}$.
The Rayleigh-Ritz procedure $\mathrm{RR}[H,s]$
extracts orthonormal Ritz vectors of $H$ associated with
the $s$ largest Ritz values.
\end{remark}

\section{Approaches and auxiliary subspaces}

In this section, we begin with exact-inverse preconditioning
$\widetilde{T}=A^{-1}$ in the BPG iteration \eqref{bpg1}
and introduce two approaches for the convergence analysis.
The first approach is a comparative analysis where the trial subspace is simplified
in order to apply estimates from \cite{KNY1987} for an abstract block iteration.
We particularly introduce some underlying auxiliary subspaces and
formulate with them the second approach. Therein proper vector iterations
are constructed for preparing the analysis for general preconditioners
based on our previous results from \cite{NEZ2014,NEZ2019}.

\subsection{Analysis via an abstract block iteration}

In the case $\widetilde{T}=A^{-1}$, we can set $N=I$
in the accompanying iteration \eqref{bpg2}. Then the trial subspace turns into
\[\spanl\{Y^{(\ell)},HY^{(\ell)}-Y^{(\ell)}\Theta_Y^{(\ell)}\}
 =\spanl\{Y^{(\ell)}\Theta_Y^{(\ell)},HY^{(\ell)}-Y^{(\ell)}\Theta_Y^{(\ell)}\}
 =\spanl\{Y^{(\ell)},HY^{(\ell)}\}\]
where the diagonal Ritz value matrix $\Theta_Y^{(\ell)}$ is invertible
due to the positive definiteness of $H$.
Therefore \eqref{bpg2} is specialized as
\begin{equation}\label{bpg2a}
 Y^{(\ell+1)} \quad\xleftarrow{\mathrm{RR}[H,s]}\quad
 \spanl\{Y^{(\ell)},HY^{(\ell)}\}.
\end{equation}
For an arbitrary linear polynomial $p_1(\cdot)$, the iteration
\begin{equation}\label{abi}
 Y^{(\ell+1)}=p_1(H)Y^{(\ell)}
\end{equation}
does not converge faster than \eqref{bpg2a}
since the Rayleigh-Ritz procedure thereof provides
the best $s$ approximate eigenvalues in the larger subspace
$\spanl\{Y^{(\ell)},HY^{(\ell)}\}$ enclosing $\spanl\{p_1(H)Y^{(\ell)}\}$.

Indeed, the iteration \eqref{bpg2a} can be regarded as
a simply restarted version of the block Lanczos method
and investigated based on the comparative analysis from \cite[Section 2]{KNY1987}
by Knyazev. We reformulate the central estimate therein as follows.

\begin{lemma}[reformulation of {\cite[(2.22)]{KNY1987}}]\label{lm1}
With the settings from Remark \ref{cas}, consider the iteration
$Y^{(\ell+1)}=f(H)Y^{(\ell)}$ for $Y^{(\ell)}\in\CC^{n \times s}$ and
a function $f(\cdot)$ satisfying $|f(\mu_1)|\ge\cdots\ge|f(\mu_s)|>0$.
If $Y^{(\ell)}$ has full rank and the $s$th largest Ritz value $\theta_s^{(\ell)}$
of $H$ in $\spanl\{Y^{(\ell)}\}$ is larger than $\mu_{s+1}$,
then $Y^{(\ell+1)}$ also has full rank. In addition, for the corresponding Ritz value
$\theta_s^{(\ell+1)}$, it holds that
\begin{equation}\label{lm1e}
 \frac{\mu_s-\theta_s^{(\ell+1)}}{\theta_s^{(\ell+1)}-\mu_{s+1}}
 \le \left(\frac{\max_{k=s+1,\ldots,n}|f(\mu_k)|}{\min_{k=1,\ldots,s}|f(\mu_k)|}\right)^2
 \frac{\mu_s-\theta_s^{(\ell)}}{\theta_s^{(\ell)}-\mu_{s+1}}.
\end{equation}
\end{lemma}

Applying Lemma \ref{lm1} to \eqref{abi} with
\[f(\mu)=p_1(\mu)=\mu-\tfrac12(\mu_{s+1}+\mu_n)\]
yields the convergence factor
\[\frac{\max_{k=s+1,\ldots,n}|f(\mu_k)|}{\min_{k=1,\ldots,s}|f(\mu_k)|}
 =\frac{|f(\mu_{s+1})|}{|f(\mu_s)|}
 =\frac{\mu_{s+1}-\mu_n}{2\mu_s-\mu_{s+1}-\mu_n}
 =\frac{\kappa_s}{2-\kappa_s}
 \quad\mbox{with}\quad
 \kappa_s=\frac{\mu_{s+1}-\mu_n}{\mu_s-\mu_n}\]
so that \eqref{lm1e} provides a single-step estimate for \eqref{bpg2a}
which can be applied recursively for multiple steps.
A direct extension to the $i$th largest Ritz value for
an arbitrary $i \le s$ does not hold in general; cf.~the numerical example
in \cite[Section 3 and Figure 1]{ZHO2018}.
In contrast, the estimate \cite[(2.20)]{KNY1987}
leads to the angle-dependent multi-step estimate
\[\frac{\mu_i-\theta_i^{(\ell)}}{\theta_i^{(\ell)}-\mu_n}
 \le \left(\frac{\kappa_i}{2-\kappa_i}\right)^{2\ell}
 \tan^2\varphi^{(0)}
 \quad\mbox{with}\quad
 \kappa_i=\frac{\mu_{s+1}-\mu_n}{\mu_i-\mu_n}\]
for \eqref{bpg2a} where $\varphi^{(0)}$ is the Euclidean angle
between the initial subspace $\spanl\{Y^{(0)}\}$
and the invariant subspace of $H$ associated with the eigenvalues
$\mu_1,\ldots,\mu_s$. Furthermore, two angle-free multi-step estimates
for \eqref{bpg2a} can be derived analogously to recent results
from \cite{NEZ2019} for the block power method $Y^{(\ell+1)}=HY^{(\ell)}$
(and swapping the indices $i$ and $j$ in the notation therein).

\begin{lemma}[based on {\cite[Theorems 2.6 and 2.8]{NEZ2019}}]\label{lm2}
With the settings from Remark \ref{cas}, consider the special form \eqref{bpg2a}
of the iteration \eqref{bpg2}. If $\theta_s^{(0)}>\mu_{s+1}$,
then it holds that
\begin{equation}\label{lm2e1}
 \frac{\mu_i-\theta_i^{(\ell)}}{\theta_i^{(\ell)}-\mu_{s+1}}
 \le \left(\frac{\kappa_i}{2-\kappa_i}\right)^{2\ell}
 \frac{\mu_i-\theta_s^{(0)}}{\theta_s^{(0)}-\mu_{s+1}}
 \quad\mbox{with}\quad
 \kappa_i=\frac{\mu_{s+1}-\mu_n}{\mu_i-\mu_n}.
\end{equation}
A more general estimate in the case $\mu_j\ge\theta_s^{(0)}>\mu_{j+1}$
with a certain index $j \ge s$ reads
\begin{equation}\label{lm2e2}
 \frac{\mu_{j-s+i}-\theta_i^{(\ell)}}{\theta_i^{(\ell)}-\mu_{j+1}}
 \le \left(\frac{\kappa_i}{2-\kappa_i}\right)^{2\ell}
 \frac{\mu_{j-s+i}-\theta_s^{(0)}}{\theta_s^{(0)}-\mu_{j+1}}
 \quad\mbox{with}\quad
 \kappa_i=\frac{\mu_{j+1}-\mu_n}{\mu_{j-s+i}-\mu_n}.
\end{equation}
\end{lemma}

For proving Lemma \ref{lm2}, we can first adapt the analysis from \cite{NEZ2019}
to the iteration \eqref{abi} and then extend the results
to the accelerated iteration \eqref{bpg2a} by using again the subspace inclusion
$\spanl\{p_1(H)Y^{(\ell)}\}\subseteq\spanl\{Y^{(\ell)},HY^{(\ell)}\}$.
An underlying proof technique is that one can select
a subspace $\widetilde{\mathcal{Y}}\subseteq\spanl\{Y^{(\ell)}\}$
such that $\widetilde{\mathcal{Y}}$, $H\widetilde{\mathcal{Y}}$ and $p_1(H)\widetilde{\mathcal{Y}}$
are simultaneously orthogonal to the eigenvectors associated with
the eigenvalues $\mu_{i+1},\ldots,\mu_s$
or $\mu_{j-s+i+1},\ldots,\mu_j$ which can be skipped in the estimates.
However, this approach cannot easily be adapted to
the iteration \eqref{bpg2} with an arbitrary $N \approx I$
except for some special cases such as
that $N$ has the same eigenvectors as $H$.
A possible way out is to construct proper
vector iterations within the subspace $\widetilde{\mathcal{Y}}+H\widetilde{\mathcal{Y}}$
and its counterpart for $N \approx I$ so that
some arguments from \cite{NEZ2014,NEZ2019} can be utilized.

\subsection{Auxiliary subspaces}

Following the proof sketch of Lemma \ref{lm2}, we introduce some auxiliary subspaces
which are still useful for analyzing general preconditioning $N \approx I$.

\begin{lemma}\label{lmaux}
With the settings from Remark \ref{cas}, let
$z_1,\ldots,z_n$ be orthonormal eigenvectors of $H$
associated with the eigenvalues $\mu_1\ge\cdots\ge\mu_n$.
By using the invariant subspaces
\[\widetilde{\mathcal{Z}}=\spanl\{z_{j-s+i+1},\ldots,z_j\}^{\perp}
 \quad\mbox{and}\quad
 \widehat{\mathcal{Z}}=\spanl\{z_{j-s+1},\ldots,z_{j-s+i}\}^{\perp}\]
for a certain index $j \ge s$ (therein the superscript 
${}^{\perp}$ denotes orthogonal complement), define for
an arbitrary subspace $\mathcal{Y}\subseteq\CC^n$
of dimension $s$ the auxiliary subspaces
\[\widetilde{\mathcal{Y}}=\mathcal{Y}\cap\widetilde{\mathcal{Z}},
 \quad\mbox{and}\quad
 \widehat{\mathcal{Y}}=\mathcal{Y}\cap\widehat{\mathcal{Z}}.\]
Then it holds that
\begin{equation}\label{lmauxe1}
 \dim\widetilde{\mathcal{Y}}\ge i,
 \quad\mbox{and}\quad
 \dim\widehat{\mathcal{Y}}\ge s{\,-\,}i.
\end{equation}
If $j=s$, and the smallest Ritz values $\widetilde{\theta}$, $\widehat{\theta}$
of $H$ in $\widetilde{\mathcal{Y}}$, $\widehat{\mathcal{Y}}$ are larger than $\mu_{s+1}$, then
\begin{equation}\label{lmauxe2}
 \dim\widetilde{\mathcal{Y}}=i,\quad
 \dim\widehat{\mathcal{Y}}=s{\,-\,}i,\quad
 \dim(\widetilde{\mathcal{Y}}\cap\widehat{\mathcal{Y}})=0,
 \quad\mbox{and}\quad
 \dim(\widetilde{\mathcal{Y}}+\widehat{\mathcal{Y}})=s.
\end{equation}
\end{lemma}
\begin{proof}
The statement \eqref{lmauxe1} follows from
\[\begin{split}
 \dim\widetilde{\mathcal{Y}}
 &=\dim\mathcal{Y}+\dim\widetilde{\mathcal{Z}}-\dim(\mathcal{Y}+\widetilde{\mathcal{Z}})
 \ge s+(n-s+i)-n=i,\\[1ex]
 \dim\widehat{\mathcal{Y}}
 &=\dim\mathcal{Y}+\dim\widehat{\mathcal{Z}}-\dim(\mathcal{Y}+\widehat{\mathcal{Z}})
 \ge s+(n-i)-n=s-i.
\end{split}\]
If $j=s$, the additional assumption on Ritz values
excludes the strict inequalities in \eqref{lmauxe1} since
\[\begin{split}
 &\dim\widetilde{\mathcal{Y}}>i
 \quad\Rightarrow\quad
 \widetilde{\theta}\le\mbox{the $(i{\,+\,}1)$th element in}\
 \{\mu_1,\ldots,\mu_i,\mu_{s+1},\ldots,\mu_n\}=\mu_{s+1},\\[1ex]
 &\dim\widehat{\mathcal{Y}}>s-i
 \quad\Rightarrow\quad
 \widehat{\theta}\le\mbox{the $(s{\,-\,}i{\,+\,}1)$th element in}\
 \{\mu_{i+1},\ldots,\mu_n\}=\mu_{s+1}.
\end{split}\]
Moreover, $\dim(\widetilde{\mathcal{Y}}\cap\widehat{\mathcal{Y}})=0$ holds
since otherwise there would exist nonzero vectors
in $\widetilde{\mathcal{Y}}\cap\widehat{\mathcal{Y}}$
and its superset $\widetilde{\mathcal{Z}}\cap\widehat{\mathcal{Z}}
=\spanl\{z_{s+1},\ldots,z_n\}$ so that
the smallest Ritz values $\widetilde{\theta}$, $\widehat{\theta}$
would be not larger than $\mu_{s+1}$. Consequently,
$\dim(\widetilde{\mathcal{Y}}+\widehat{\mathcal{Y}})
=\dim\widetilde{\mathcal{Y}}+\dim\widehat{\mathcal{Y}}
-\dim(\widetilde{\mathcal{Y}}\cap\widehat{\mathcal{Y}})=s$.
\end{proof}

Based on the statement \eqref{lmauxe1} and the Courant-Fischer principles,
the convergence of the $i$th Ritz value produced by \eqref{bpg2a}
is not slower than that of the $i$th Ritz value by the iteration
\begin{equation}\label{bpg2a1}
 \widetilde{Y}^{(\ell+1)} \quad\xleftarrow{\mathrm{RR}[H,\widetilde{i}]}\quad
 \spanl\{\widetilde{Y}^{(\ell)},H\widetilde{Y}^{(\ell)}\}
 \quad\mbox{with}\quad
 \spanl\{\widetilde{Y}^{(0)}\}=\spanl\{Y^{(0)}\}\cap\widetilde{\mathcal{Z}}
\end{equation}
and $\widetilde{i}=\dim\spanl\{\widetilde{Y}^{(0)}\}\ge i$.
Evidently, each iterate of \eqref{bpg2a1} is contained columnwise
in $\widetilde{\mathcal{Z}}$ so that the estimate \eqref{lm2e2} can be derived
by modifying Lemma \ref{lm1} restricted to $\widetilde{\mathcal{Z}}$.
The statement \eqref{lmauxe2} refers to the final phase of \eqref{bpg2a}
and the estimate \eqref{lm2e1} under the assumption $\theta_s^{(0)}>\mu_{s+1}$
(then the corresponding $\widetilde{\theta}$ and $\widehat{\theta}$
are also larger than $\mu_{s+1}$). Therein \eqref{bpg2a} can be split into
two partial iterations with respect to $\widetilde{\mathcal{Z}}$ and $\widehat{\mathcal{Z}}$.
Moreover, we can inductively adapt \eqref{lmauxe2}
to the respective subspace iterates; see Lemma \ref{lmaux1}.

\subsection{Analysis via vector iterations}

For analyzing the iteration \eqref{bpg2} with general preconditioners,
the direct generalization
\[\widetilde{Y}^{(\ell+1)} \quad\xleftarrow{\mathrm{RR}[H,\widetilde{i}]}\quad
 \spanl\{\widetilde{Y}^{(\ell)},N(H\widetilde{Y}^{(\ell)}-\widetilde{Y}^{(\ell)}\Theta_{\widetilde{Y}}^{(\ell)})\}\]
of \eqref{bpg2a1} is somewhat problematic since
$\spanl\{N(H\widetilde{Y}^{(\ell)}-\widetilde{Y}^{(\ell)}\Theta_{\widetilde{Y}}^{(\ell)})\}$
for $N \approx I$ is not necessarily a subset of $\widetilde{\mathcal{Z}}$.
Instead, following our previous results from \cite{NEZ2014,NEZ2019},
we reformulate the trial subspace of \eqref{bpg2} as
\[\spanl\{Y^{(\ell)},N(HY^{(\ell)}-Y^{(\ell)}\Theta_Y^{(\ell)})\}
=\spanl\{Y^{(\ell)},Y^{(\ell)}\Theta_Y^{(\ell)}+N(HY^{(\ell)}-Y^{(\ell)}\Theta_Y^{(\ell)})\}\]
and consider a stepwise mixture of two partial iterations concerning
$\widetilde{\mathcal{Z}}$ and $\widehat{\mathcal{Z}}$, namely,
\begin{equation}\label{bpg2a2}
\begin{split}
 &\widetilde{\mathcal{Y}}^{(\ell)}=\spanl\{Y^{(\ell)}\}\cap\widetilde{\mathcal{Z}},\qquad
 \widehat{\mathcal{Y}}^{(\ell)}=\spanl\{Y^{(\ell)}\}\cap\widehat{\mathcal{Z}},\\[1ex]
 &U^{(\ell)}=Y^{(\ell)}\Theta_Y^{(\ell)}+N(HY^{(\ell)}-Y^{(\ell)}\Theta_Y^{(\ell)}),\\[1ex]
 &\,\widetilde{\mathcal{U}}^{(\ell)}=\spanl\{U^{(\ell)}\}\cap\widetilde{\mathcal{Z}},\qquad
 \,\widehat{\mathcal{U}}^{(\ell)}=\spanl\{U^{(\ell)}\}\cap\widehat{\mathcal{Z}},\\[1ex]
 &\widetilde{Y}^{(\ell+\frac12)} \quad\xleftarrow{\mathrm{RR}[H,\widetilde{i}]}\quad
 \widetilde{\mathcal{Y}}^{(\ell)}+\,\widetilde{\mathcal{U}}^{(\ell)},\qquad
 \widehat{Y}^{(\ell+\frac12)} \quad\xleftarrow{\mathrm{RR}[H,\widehat{i}]}\quad
 \widehat{\mathcal{Y}}^{(\ell)}+\,\widehat{\mathcal{U}}^{(\ell)},\\[1ex]
 &Y^{(\ell+1)} \quad\xleftarrow{\mathrm{RR}[H,\widetilde{s}]}\quad
 \spanl\{\widetilde{Y}^{(\ell+\frac12)},\widehat{Y}^{(\ell+\frac12)}\}
\end{split}
\end{equation}
with $\widetilde{i}=\dim\widetilde{\mathcal{Y}}^{(\ell)}$, $\widehat{i}=\dim\widehat{\mathcal{Y}}^{(\ell)}$
and $\widetilde{s}=\dim\spanl\{\widetilde{Y}^{(\ell+\frac12)},\widehat{Y}^{(\ell+\frac12)}\}$
for the current step index $\ell$. The matrix $U^{(\ell)}$ coincides with
$HY^{(\ell)}$ for $N=I$, and $\spanl\{U^{(\ell)}\}$ corresponds to 
the trial subspace of a BPG iteration with fixed step sizes
for which some cluster robust estimates have been derived in \cite{NEZ2019}.
The trial subspace of \eqref{bpg2}, i.e., $\spanl\{Y^{(\ell)},U^{(\ell)}\}$,
is split into $\widetilde{\mathcal{Y}}^{(\ell)}+\,\widetilde{\mathcal{U}}^{(\ell)}$
and $\widehat{\mathcal{Y}}^{(\ell)}+\,\widehat{\mathcal{U}}^{(\ell)}$
which are subsets of $\widetilde{\mathcal{Z}}$ and $\widehat{\mathcal{Z}}$, respectively.
Therein two partial Rayleigh-Ritz approximations are determined
and additionally refined together for extracting the next iterate.
In comparison to the direct Rayleigh-Ritz approximation 
in the larger subspace $\spanl\{Y^{(\ell)},U^{(\ell)}\}$ in \eqref{bpg2},
the update by \eqref{bpg2a2} leads to less improvement in
Ritz values according to the Courant-Fischer principles.
Thus investigating \eqref{bpg2a2} can provide suitable Ritz value estimates for \eqref{bpg2}.
The next task in this approach is to construct proper vector iterations within \eqref{bpg2a2}
as well as an alternative quality parameter for $N \approx I$.

We first discuss the dimensions of auxiliary subspaces
in \eqref{bpg2a2} for $j{\,=\,}s$
concerning the final phase of the iteration \eqref{bpg2}.

\begin{lemma}\label{lmaux1}
With the settings from Remark \ref{cas},
consider the $\ell$th step of \eqref{bpg2a2} with the invariant subspaces
$\widetilde{\mathcal{Z}}$ and $\widehat{\mathcal{Z}}$ for $j{\,=\,}s$ from Lemma \ref{lmaux}.
If $\dim\widetilde{\mathcal{Y}}^{(\ell)}=i$, $\dim\widehat{\mathcal{Y}}^{(\ell)}=s{\,-\,}i$,
and the smallest Ritz values $\widetilde{\theta}$, $\widehat{\theta}$
of $H$ in $\widetilde{\mathcal{Y}}^{(\ell)}$, $\widehat{\mathcal{Y}}^{(\ell)}$
are larger than $\mu_{s+1}$, then it holds for the subspaces
$\widetilde{\mathcal{Y}}'=\spanl\{\widetilde{Y}^{(\ell+\frac12)}\}$ and
$\widehat{\mathcal{Y}}'=\spanl\{\widehat{Y}^{(\ell+\frac12)}\}$, that
$\dim\widetilde{\mathcal{Y}}'=i$, $\dim\widehat{\mathcal{Y}}'=s{\,-\,}i$,
and the smallest Ritz values $\widetilde{\theta}'$, $\widehat{\theta}'$
of $H$ in $\widetilde{\mathcal{Y}}'$, $\widehat{\mathcal{Y}}'$
are also larger than $\mu_{s+1}$. Moreover,
\begin{equation}\label{lmaux1e}
 \dim(\widetilde{\mathcal{Y}}'+\widehat{\mathcal{Y}}')=s,\quad
 \widetilde{\mathcal{Y}}^{(\ell+1)}=\widetilde{\mathcal{Y}}',
 \quad\mbox{and}\quad
 \widehat{\mathcal{Y}}^{(\ell+1)}=\widehat{\mathcal{Y}}'.
\end{equation} 
\end{lemma}
\begin{proof}
The given assumption leads to $\widetilde{i}=i$ and $\widehat{i}=s{\,-\,}i$ so that
the partial Rayleigh-Ritz approximations produce
$\widetilde{\mathcal{Y}}'$ of dimension $i$ and
$\widehat{\mathcal{Y}}'$ of dimension $s{\,-\,}i$.
In addition, the smallest Ritz values $\widetilde{\theta}'$, $\widehat{\theta}'$
of $H$ in $\widetilde{\mathcal{Y}}'$, $\widehat{\mathcal{Y}}'$
improve $\widetilde{\theta}$, $\widehat{\theta}$, namely,
\[\widetilde{\theta}'=\theta_i(\widetilde{\mathcal{Y}}^{(\ell)}+\,\widetilde{\mathcal{U}}^{(\ell)})
 \ge\theta_i(\widetilde{\mathcal{Y}}^{(\ell)})=\widetilde{\theta},
 \quad\mbox{and}\quad
 \widehat{\theta}'=\theta_{s-i}(\widehat{\mathcal{Y}}^{(\ell)}+\,\widehat{\mathcal{U}}^{(\ell)})
 \ge\theta_{s-i}(\widehat{\mathcal{Y}}^{(\ell)})=\widehat{\theta}.\]
Thus $\widetilde{\theta}'$ and $\widehat{\theta}'$ are also larger than $\mu_{s+1}$.
Subsequently, the property $\dim(\widetilde{\mathcal{Y}}'+\widehat{\mathcal{Y}}')=s$
can be shown analogously to the last equality in \eqref{lmauxe2}.
Then $\widetilde{s}=s$, and
\[\widetilde{\mathcal{Y}}^{(\ell+1)}=\spanl\{Y^{(\ell+1)}\}\cap\widetilde{\mathcal{Z}}
 =(\widetilde{\mathcal{Y}}'\cap\widetilde{\mathcal{Z}})+(\widehat{\mathcal{Y}}'\cap\widetilde{\mathcal{Z}})
 =\widetilde{\mathcal{Y}}'\]
holds according to
\[\begin{split}
 &\widetilde{\mathcal{Y}}'\subseteq(\widetilde{\mathcal{Y}}^{(\ell)}{\,+\,}\,\widetilde{\mathcal{U}}^{(\ell)})
 \subseteq\widetilde{\mathcal{Z}}\quad\Rightarrow\quad
 \widetilde{\mathcal{Y}}'\cap\widetilde{\mathcal{Z}}=\widetilde{\mathcal{Y}}',\\[1ex]
 &\widehat{\mathcal{Y}}'\subseteq(\widehat{\mathcal{Y}}^{(\ell)}{\,+\,}\,\widehat{\mathcal{U}}^{(\ell)})
 \subseteq\widehat{\mathcal{Z}}\quad\Rightarrow\quad
 \widehat{\mathcal{Y}}'\cap\widehat{\mathcal{Z}}=\widehat{\mathcal{Y}}'\\[1ex]
 &\hspace{1cm}\quad\Rightarrow\quad
 \widehat{\mathcal{Y}}'\cap\widetilde{\mathcal{Z}}
 =\widehat{\mathcal{Y}}'\cap\widehat{\mathcal{Z}}\cap\widetilde{\mathcal{Z}}
 =\widehat{\mathcal{Y}}'\cap\spanl\{z_{s+1},\ldots,z_n\}=\{0\}
\end{split}\]
where the last equality is ensured by $\widehat{\theta}'>\mu_{s+1}$.
The verification of \,$\widehat{\mathcal{Y}}^{(\ell+1)}=\widehat{\mathcal{Y}}'$ is analogous.
\end{proof}

Lemma \ref{lmaux1} enables an inductive proof of the following properties
of \eqref{bpg2a2} under a natural assumption on the initial subspace.

\begin{lemma}\label{lmaux2}
With the settings from Remark \ref{cas},
consider the iteration \eqref{bpg2a2} with the invariant subspaces
$\widetilde{\mathcal{Z}}$ and $\widehat{\mathcal{Z}}$ for $j{\,=\,}s$ from Lemma \ref{lmaux}
and $\dim\spanl\{Y^{(0)}\}=s$. If the smallest ($s$th largest) Ritz value
$\theta_s^{(0)}$ of $H$ in $\spanl\{Y^{(0)}\}$ is larger than $\mu_{s+1}$,
then it holds for each $\ell$, that $\dim\widetilde{\mathcal{Y}}^{(\ell)}=i$,
$\dim\widehat{\mathcal{Y}}^{(\ell)}=s{\,-\,}i$, and $\dim\spanl\{Y^{(\ell+1)}\}=s$.
The partial Rayleigh-Ritz approximations in the $\ell$th step
actually produce the subspaces
$\widetilde{\mathcal{Y}}^{(\ell+1)}$ and $\widehat{\mathcal{Y}}^{(\ell+1)}$.
\end{lemma}
\begin{proof}
Applying Lemma \ref{lmaux} to $\mathcal{Y}=\spanl\{Y^{(0)}\}$
implies $\dim\widetilde{\mathcal{Y}}^{(0)}=i$ and $\dim\widehat{\mathcal{Y}}^{(0)}=s{\,-\,}i$
by the first two equalities in \eqref{lmauxe2}
whose assumption is verified by the fact
that the smallest Ritz values $\widetilde{\theta}$ and $\widehat{\theta}$
of $H$ in the subsets $\widetilde{\mathcal{Y}}^{(0)}$ and $\widehat{\mathcal{Y}}^{(0)}$
of $\spanl\{Y^{(0)}\}$ are at least $\theta_s^{(0)}$
and thus larger than $\mu_{s+1}$. Therefore Lemma \ref{lmaux1}
is already applicable to $\ell=0$. Moreover, the statements for $\ell$
in Lemma \ref{lmaux1} immediately verify the assumption for $\ell{\,+\,}1$.
Recursively applying Lemma \ref{lmaux1} completes the proof.
\end{proof}

Lemma \ref{lmaux2} motivates an approach for estimating the convergence rate
of the $i$th Ritz value in the final phase of the iteration \eqref{bpg2}
by observing the partial subspace iterate $\widetilde{\mathcal{Y}}^{(\ell)}$ in \eqref{bpg2a2}.
The other partial subspace iterate $\widehat{\mathcal{Y}}^{(\ell)}$ plays
an important role in the background for ensuring $\dim\spanl\{Y^{(\ell+1)}\}=s$.
Extending Lemma \ref{lmaux2} to the more general case $j \ge s$ requires certain
assumptions on the initial subspace $\spanl\{Y^{(0)}\}$ which are
much more technical than the natural assumption $\theta_s^{(0)}>\mu_{s+1}$.
It is remarkable that opposite properties such as
\,$\dim\widetilde{\mathcal{Y}}^{(\ell)}>i$\,
rarely occur in numerical tests with randomly generated initial guesses.
Therefore we simply use an empirical assumption for analyzing \eqref{bpg2a2}
in the case $j \ge s$.

\begin{lemma}\label{lmaux3}
With the settings from Remark \ref{cas},
consider the iteration \eqref{bpg2a2} with the invariant subspaces
$\widetilde{\mathcal{Z}}$ and $\widehat{\mathcal{Z}}$ for $j{\,\ge\,}s$ from Lemma \ref{lmaux}.
Assume for each $\ell$ that $\dim\spanl\{Y^{(\ell)}\}=s$,
$\dim\widetilde{\mathcal{Y}}^{(\ell)}=i$, and $\dim\widehat{\mathcal{Y}}^{(\ell)}=s{\,-\,}i$.
Then the partial Rayleigh-Ritz approximations in the $\ell$th step
actually produce the subspaces
$\widetilde{\mathcal{Y}}^{(\ell+1)}$ and $\widehat{\mathcal{Y}}^{(\ell+1)}$.
\end{lemma}
\begin{proof}
The statement cannot be proved by directly applying Lemma \ref{lmaux1}
due to the dependence on Ritz values. Instead,
\[\widetilde{\mathcal{Y}}^{(\ell+1)}=(\spanl\{Y^{(\ell+1)}\}\cap\widetilde{\mathcal{Z}})
 \ \supseteq \ (\widetilde{\mathcal{Y}}'\cap\widetilde{\mathcal{Z}})=\widetilde{\mathcal{Y}}'
 \quad\Rightarrow\quad \widetilde{\mathcal{Y}}^{(\ell+1)}=\widetilde{\mathcal{Y}}'\]
holds since $\dim\widetilde{\mathcal{Y}}'=i$ and
$\dim\widetilde{\mathcal{Y}}^{(\ell+1)}=i$ (by adapting the assumption to $\ell{\,+\,}1$).
The equality $\widehat{\mathcal{Y}}^{(\ell+1)}=\widehat{\mathcal{Y}}'$ holds analogously.
\end{proof}

Now we can focus on the first partial iteration in \eqref{bpg2a2}
and define certain vector iterations for characterizing the $i$th Ritz value.

\begin{theorem}\label{thm1}
With the settings from Remark \ref{cas},
consider the iteration \eqref{bpg2a2} under the assumption
from Lemma \ref{lmaux2} or Lemma \ref{lmaux3}, 
and denote by $\widetilde{\theta}'$ the $i$th largest Ritz value
of $H$ in $\widetilde{\mathcal{Y}}^{(\ell)}{\,+\,}\,\widetilde{\mathcal{U}}^{(\ell)}$.
Then the following statements hold:
\begin{itemize}
\setlength\itemsep{1ex}
\item[(a)] For each $\ell$, the matrix $U^{(\ell)}$ has full rank,
and $\dim\,\widetilde{\mathcal{U}}^{(\ell)} \ge i$.
\item[(b)] In the special case $N=I$, the subspace
$\,\widetilde{\mathcal{U}}^{(\ell)}$ coincides with $H\widetilde{\mathcal{Y}}^{(\ell)}$,
and there exists a nonzero vector $\widetilde{y}\in\widetilde{\mathcal{Y}}^{(\ell)}$
for which the largest Ritz value $\widetilde{\theta}^{\diamond}$ of $H$ in $\spanl\{\widetilde{y},H\widetilde{y}\}$
does not exceed $\widetilde{\theta}'$.
\item[(c)] In the general case $N \approx I$,
consider an orthonormal matrix $\widetilde{U}\in\CC^{n \times i}$
with $\spanl\{\widetilde{U}\}\subseteq\widetilde{\mathcal{U}}^{(\ell)}$,
and an orthonormal basis matrix $\widetilde{Y}$ of $\widetilde{\mathcal{Y}}^{(\ell)}$.
Let $\mu(\cdot)$ be the Rayleigh quotient with respect to $H$.
If the matrix $\widetilde{R}=H\widetilde{Y}-\widetilde{Y}\widetilde{Y}^*H\widetilde{Y}$
has full rank, and
\begin{equation}\label{prec2a}
 \|(H\widetilde{Y}-\widetilde{U}\widetilde{U}^*H\widetilde{Y})(\widetilde{R}^*\widetilde{R})^{-1/2}\|_2
 \le\widetilde{\gamma}<1,
\end{equation}
then there exist nonzero vectors $\widetilde{y}\in\spanl\{\widetilde{Y}\}$
and $\widetilde{u}\in\spanl\{\widetilde{U}\}$ such that
\begin{equation}\label{prec2b}
 \|H\widetilde{y}-\widetilde{u}\|_2\le\widetilde{\gamma}\|H\widetilde{y}-\mu(\widetilde{y})\widetilde{y}\|_2,
\end{equation}
and the largest Ritz value $\widetilde{\theta}^{\diamond}$ of $H$ in $\spanl\{\widetilde{y},\widetilde{u}\}$
does not exceed $\widetilde{\theta}'$.
\end{itemize}
\end{theorem}
\begin{proof}
(a) According to Lemma \ref{lmaux2} or Lemma \ref{lmaux3},
we get $\dim\spanl\{Y^{(\ell)}\}=s$ for each $\ell$.
The corresponding $U^{(\ell)}$ can be represented by
\[\begin{split}
 &U^{(\ell)}=Y^{(\ell)}\Theta_Y^{(\ell)}+N(HY^{(\ell)}-Y^{(\ell)}\Theta_Y^{(\ell)})
 =Y'\Theta\\[1ex]
 \quad\mbox{with}\quad
 &Y=Y^{(\ell)}, \quad \Theta=\Theta_Y^{(\ell)}
 \quad\mbox{and}\quad
 Y'=Y-N(Y-HY\Theta^{-1})
\end{split}\]
for matching the notation in \cite[Lemma 3.1]{NEZ2019}
where a BPG iteration with fixed step sizes is analyzed,
and $Y'$ can be shown to have full rank.
Then $U^{(\ell)}=Y'\Theta$ also has full rank
since the diagonal Ritz value matrix $\Theta=\Theta_Y^{(\ell)}$
is invertible due to the positive definiteness of $H$.
Subsequently, $\dim\,\widetilde{\mathcal{U}}^{(\ell)} \ge i$
can be shown analogously to \eqref{lmauxe1} in Lemma \ref{lmaux}.

(b) For $N=I$, the matrix $U^{(\ell)}$ becomes $HY^{(\ell)}$ so that
\[\,\widetilde{\mathcal{U}}^{(\ell)}=\spanl\{U^{(\ell)}\}\cap\widetilde{\mathcal{Z}}
 =\spanl\{HY^{(\ell)}\}\cap(H\widetilde{\mathcal{Z}})
 =H(\spanl\{Y^{(\ell)}\}\cap\widetilde{\mathcal{Z}})
 =H\widetilde{\mathcal{Y}}^{(\ell)}\]
(where $\widetilde{\mathcal{Z}}=H\widetilde{\mathcal{Z}}$
is ensured by the positive definiteness of $H$).
Following the property $\dim\widetilde{\mathcal{Y}}^{(\ell)}=i$
from Lemma \ref{lmaux2} or Lemma \ref{lmaux3},
we use an arbitrary basis matrix $\widetilde{Y}\in\CC^{n \times i}$
of $\widetilde{\mathcal{Y}}^{(\ell)}$ so that the subspace
$\,\widetilde{\mathcal{U}}'=\widetilde{\mathcal{Y}}^{(\ell)}{\,+\,}\,\widetilde{\mathcal{U}}^{(\ell)}$
can be represented by $\spanl\{\widetilde{Y},H\widetilde{Y}\}$.
We denote by $\bar{i}$ the dimension of $\,\widetilde{\mathcal{U}}'$,
and by $V\in\CC^{n \times \bar{i}}$ a basis matrix of $\,\widetilde{\mathcal{U}}'$
whose columns $v_1,\ldots,v_{\,\bar{i}}$ are orthonormal Ritz vectors
associated with the Ritz values
$\varphi_1\ge\cdots\ge\varphi_{\,\bar{i}}$ of $H$ in $\,\widetilde{\mathcal{U}}'$.
Then we get the orthogonal projector $P=VV^*$ on $\,\widetilde{\mathcal{U}}'$
and the diagonal Ritz value matrix $V^*HV=\diag(\varphi_1,\ldots,\varphi_{\,\bar{i}})$.
Moreover, the $i$th largest Ritz value $\widetilde{\theta}'$
of $H$ in $\,\widetilde{\mathcal{U}}'$ is the largest Ritz value of $H$
in $\,\widetilde{\mathcal{U}}^{\diamond}=\spanl\{v_i,\ldots,v_{\,\bar{i}}\}$.
Based on the dimension comparison
\[\dim(\widetilde{\mathcal{Y}}^{(\ell)}\cap\,\widetilde{\mathcal{U}}^{\diamond})
 =\dim\widetilde{\mathcal{Y}}^{(\ell)}+\dim\,\widetilde{\mathcal{U}}^{\diamond}
 -\dim(\widetilde{\mathcal{Y}}^{(\ell)}+\,\widetilde{\mathcal{U}}^{\diamond})
 \ge i+(\bar{i}-i+1)-\bar{i}=1,\]
we can select a nonzero vector $\widetilde{y}$
from $\widetilde{\mathcal{Y}}^{(\ell)}\cap\,\widetilde{\mathcal{U}}^{\diamond}$.
According to $\widetilde{y}\in\widetilde{\mathcal{Y}}^{(\ell)}$ and
$H\widetilde{y}\in H\widetilde{\mathcal{Y}}^{(\ell)}=\,\widetilde{\mathcal{U}}^{(\ell)}$,
the vectors $\widetilde{y}$ and $H\widetilde{y}$ are contained in $\,\widetilde{\mathcal{U}}'$ so that
\[H\widetilde{y}=P(H\widetilde{y})=PH(P\widetilde{y})
 =VV^*HVV^*\widetilde{y}=V\,\diag(\varphi_1,\ldots,\varphi_{\,\bar{i}})\,V^*\widetilde{y}.\]
In addition, $\widetilde{y}\in\,\widetilde{\mathcal{U}}^{\diamond}$ and the orthogonality
between the columns of $V$ ensure that
the first $i{\,-\,}1$ entries of $V^*\widetilde{y}$ are equal to zero.
This property is preserved in the vector
$\,\diag(\varphi_1,\ldots,\varphi_{\,\bar{i}})\,V^*\widetilde{y}\,$
so that $H\widetilde{y}$ belongs to $\,\widetilde{\mathcal{U}}^{\diamond}$.
Therefore $\spanl\{\widetilde{y},H\widetilde{y}\}$ is a subset
of $\,\widetilde{\mathcal{U}}^{\diamond}$, and the largest Ritz value $\widetilde{\theta}^{\diamond}$
of $H$ in $\spanl\{\widetilde{y},H\widetilde{y}\}$ is bounded from above
by $\widetilde{\theta}'$ which is the largest Ritz value 
of $H$ in $\,\widetilde{\mathcal{U}}^{\diamond}$.

(c) The existence of $\widetilde{U}$ follows from (a).
Vectors $\widetilde{y}$ and $\widetilde{u}$ satisfying \eqref{prec2b}
can be constructed by using an arbitrary nonzero vector $c\in\CC^{i}$, namely,
\eqref{prec2a} ensures
\[\|(H\widetilde{Y}-\widetilde{U}\widetilde{U}^*H\widetilde{Y})(\widetilde{R}^*\widetilde{R})^{-1/2}c\|_2
 \le\widetilde{\gamma}\|c\|_2\]
so that
\[\|(H\widetilde{Y}-\widetilde{U}\widetilde{U}^*H\widetilde{Y})e\|_2
 \le\widetilde{\gamma}\|(\widetilde{R}^*\widetilde{R})^{1/2}e\|_2
 \quad\mbox{for}\quad
 e=(\widetilde{R}^*\widetilde{R})^{-1/2}c.\]
Subsequently, by using
$\|(\widetilde{R}^*\widetilde{R})^{1/2}e\|_2=\sqrt{e^*\widetilde{R}^*\widetilde{R}e\,}
=\|\widetilde{R}e\|_2$ and the definition of $\widetilde{R}$, we get
\[\|H\widetilde{Y}e-\widetilde{U}\widetilde{U}^*H\widetilde{Y}e\|_2
 \le\widetilde{\gamma}\|H\widetilde{Y}e-\widetilde{Y}\widetilde{Y}^*H\widetilde{Y}e\|_2
 \le\widetilde{\gamma}\|H\widetilde{Y}e-\widetilde{Y}e\,\mu(\widetilde{Y}e)\|_2\]
where the second inequality uses the fact that $\widetilde{Y}\widetilde{Y}^*H\widetilde{Y}e$
is the orthogonal projection of $H\widetilde{Y}e$ on $\spanl\{\widetilde{Y}\}$.
Thus \eqref{prec2b} is fulfilled by $\widetilde{y}=\widetilde{Y}e$
and $\widetilde{u}=\widetilde{U}(\widetilde{U}^*H\widetilde{Y}e)$.
Specific $\widetilde{y}$ and $\widetilde{u}$ possessing the additional property
can be constructed analogously to the proof of \cite[Theorem 3.2]{NEZ2014}
(using Sion's-minimax theorem).
Therein the largest Ritz value $\widetilde{\theta}^{\diamond}$
of $H$ in $\spanl\{\widetilde{y},\widetilde{u}\}$ does not exceed
the $i$th largest Ritz value $\widetilde{\theta}^{\circ}$
of $H$ in $\spanl\{\widetilde{Y},\widetilde{U}\}$.
Consequently, we get $\widetilde{\theta}^{\diamond}\le\widetilde{\theta}^{\circ}\le\widetilde{\theta}'$ according to
$\spanl\{\widetilde{Y},\widetilde{U}\}
\subseteq(\widetilde{\mathcal{Y}}^{(\ell)}{\,+\,}\,\widetilde{\mathcal{U}}^{(\ell)})$
and the Courant-Fischer principles.
\end{proof}

The statement (b) in Theorem \ref{thm1} suggests the vector iteration
\begin{equation}\label{vit1}
 \widetilde{y}^{\diamond} \quad\xleftarrow{\mathrm{RR}[H,1]}\quad
 \spanl\{\widetilde{y},H\widetilde{y}\}
\end{equation} 
for deriving an intermediate estimate. Since $\widetilde{y}$ and $H\widetilde{y}$
are contained in the invariant subspace $\widetilde{\mathcal{Z}}$,
we adapt an estimate from \cite[Theorem 4.1]{NOZ2011}
for vectorial gradient iterations as follows.

\begin{lemma}\label{lmvit1}
With the settings from Remark \ref{cas}, consider the iteration \eqref{vit1},
and let $\mu(\cdot)$ be the Rayleigh quotient with respect to $H$.
If $\widetilde{y}$ belongs to the invariant subspace $\widetilde{\mathcal{Z}}$
defined in Lemma \ref{lmaux}, and $\mu(\widetilde{y})$ is located
in the eigenvalue interval $(\mu_{j+1},\mu_{j-s+i}]$, then it holds that
\begin{equation}\label{vit1e}
 \frac{\mu_{j-s+i}-\mu(\widetilde{y}^{\diamond})}{\mu(\widetilde{y}^{\diamond})-\mu_{j+1}}
 \le \left(\frac{\kappa_i}{2-\kappa_i}\right)^2
 \frac{\mu_{j-s+i}-\mu(\widetilde{y})}{\mu(\widetilde{y})-\mu_{j+1}}
 \quad\mbox{with}\quad
 \kappa_i=\frac{\mu_{j+1}-\mu_n}{\mu_{j-s+i}-\mu_n}.
\end{equation} 
\end{lemma}
\begin{proof}
The iteration \eqref{vit1} is equivalent to
\[\widetilde{Z}^*\widetilde{y}^{\diamond} \quad\xleftarrow{\mathrm{RR}[\widetilde{Z}^*H\widetilde{Z},1]}\quad
 \spanl\{\widetilde{Z}^*\widetilde{y},\,(\widetilde{Z}^*H\widetilde{Z})\widetilde{Z}^*\widetilde{y}\}\]
with the orthonormal basis matrix $\widetilde{Z}=[z_1,\ldots,z_{j-s+i},z_{j+1},\ldots,z_n]$
of $\widetilde{\mathcal{Z}}$. Then \eqref{vit1e} is achieved by adapting
\cite[Theorem 4.1]{NOZ2011} to the matrix $\widetilde{Z}^*H\widetilde{Z}$
and the corresponding Rayleigh quotient $\widetilde{\mu}(\cdot)$
together with simple reformulations based on
\[\widetilde{\mu}(\widetilde{Z}^*w)
 =\frac{(\widetilde{Z}^*w)^*(\widetilde{Z}^*H\widetilde{Z})(\widetilde{Z}^*w)}{(\widetilde{Z}^*w)^*(\widetilde{Z}^*w)}
 =\frac{(\widetilde{Z}\widetilde{Z}^*w)^*H(\widetilde{Z}\widetilde{Z}^*w)}{w^*(\widetilde{Z}\widetilde{Z}^*w)}
 =\frac{w^*Hw}{w^*w}=\mu(w)\]
for arbitrary nonzero vectors $w$ from $\widetilde{\mathcal{Z}}$.
\end{proof}

A similar intermediate estimate for $N \approx I$ can be derived within the iteration
\begin{equation}\label{vit2}
 \widetilde{y}^{\diamond} \quad\xleftarrow{\mathrm{RR}[H,1]}\quad
 \spanl\{\widetilde{y},\widetilde{u}\}
\end{equation} 
suggested by the statement (c) in Theorem \ref{thm1}.
The derivation is essentially based on \cite{NEY2012}.

\begin{lemma}\label{lmvit2}
With the settings from Remark \ref{cas}, consider the iteration \eqref{vit2},
and let $\mu(\cdot)$ be the Rayleigh quotient with respect to $H$.
If $\widetilde{y}$ and $\widetilde{u}$ belong to the invariant subspace $\widetilde{\mathcal{Z}}$
defined in Lemma \ref{lmaux} and satisfy the condition \eqref{prec2b}
with a certain $\widetilde{\gamma}\in[0,1)$, and $\mu(\widetilde{y})$ is located
in the eigenvalue interval $(\mu_{j+1},\mu_{j-s+i}]$, then it holds that
\begin{equation}\label{vit2e}
 \frac{\mu_{j-s+i}-\mu(\widetilde{y}^{\diamond})}{\mu(\widetilde{y}^{\diamond})-\mu_{j+1}}
 \le \left(\frac{\kappa_i+\widetilde{\gamma}(2-\kappa_i)}
 {(2-\kappa_i)+\widetilde{\gamma}\kappa_i}\right)^2
 \frac{\mu_{j-s+i}-\mu(\widetilde{y})}{\mu(\widetilde{y})-\mu_{j+1}}
 \quad\mbox{with}\quad
 \kappa_i=\frac{\mu_{j+1}-\mu_n}{\mu_{j-s+i}-\mu_n}.
\end{equation} 
\end{lemma}
\begin{proof}
As in the proof of Lemma \ref{lmvit1}, we define the matrix $\widetilde{H}=\widetilde{Z}^*H\widetilde{Z}$
and the corresponding Rayleigh quotient $\widetilde{\mu}(\cdot)$ so that \eqref{vit2}
is equivalent to
\[\check{y}^{\diamond} \quad\xleftarrow{\mathrm{RR}[\widetilde{H},1]}\quad
 \spanl\{\check{y},\check{u}\}
 \quad\mbox{with}\quad
 \check{y}^{\diamond}=\widetilde{Z}^*\widetilde{y}^{\diamond},\quad
 \check{y}=\widetilde{Z}^*\widetilde{y}
 \quad\mbox{and}\quad
 \check{u}=\widetilde{Z}^*\widetilde{u}.\]
The condition \eqref{prec2b} can be reformulated as
\[\|\widetilde{H}\check{y}-\check{u}\|_2
 \le\widetilde{\gamma}\|\widetilde{H}\check{y}-\widetilde{\mu}(\check{y})\check{y}\|_2\]
since $\|w\|_2=\|\widetilde{Z}\widetilde{Z}^*w\|_2=\|\widetilde{Z}^*w\|_2$ holds
for arbitrary $w\in\widetilde{\mathcal{Z}}$. Thus $\check{u}$ belongs to
a ball $\mathcal{B}_{\widetilde{\gamma},\check{y}}$
centred at $\widetilde{H}\check{y}$ with the radius
$\widetilde{\gamma}\|\widetilde{H}\check{y}-\widetilde{\mu}(\check{y})\check{y}\|_2$.
Then the trial subspace $\spanl\{\check{y},\check{u}\}$
is characterized by a cone as in \cite[Section 2]{NEY2012}
so that the geometric analysis therefrom is applicable.
Adapting \cite[Theorem 2.2]{NEY2012} yields \eqref{vit2e}.
\end{proof}

\section{Main results}

The analysis of the auxiliary iteration \eqref{bpg2a2}
via vector iterations from Subsection 3.3
results in multi-step estimates for \eqref{bpg2} in Theorem \ref{thm2}
and corresponding estimates for the BPG iteration \eqref{bpg1} in Theorem \ref{thm3}.
The results are formulated for general preconditioning
and contain estimates for $N=I$ as special forms.

\begin{theorem}\label{thm2}
With the settings from Remark \ref{cas} concerning the iteration \eqref{bpg2},
let $z_1,\ldots,z_n$ be orthonormal eigenvectors of $H$
associated with the eigenvalues $\mu_1\ge\cdots\ge\mu_n$.
Then the following statements hold:
\begin{itemize}
\setlength\itemsep{1ex}
\item[(a)] The Ritz values produced by \eqref{bpg2} fulfill $\theta_i^{(\ell+1)}\ge\theta_i^{(\ell)}$
for each $\ell$ and $i\in\{1,\ldots,s\}$. If there are no eigenvectors in $\spanl\{Y^{(\ell)}\}$,
then $\theta_i^{(\ell+1)}>\theta_i^{(\ell)}$.
\item[(b)] If $\theta_s^{(0)}>\mu_{s+1}$, consider the auxiliary iteration \eqref{bpg2a2}
using the same initial subspace $\spanl\{Y^{(0)}\}$ together with the invariant subspaces
$\widetilde{\mathcal{Z}}=\spanl\{z_{i+1},\ldots,z_s\}^{\perp}$
and $\widehat{\mathcal{Z}}=\spanl\{z_1,\ldots,z_i\}^{\perp}$.
Then $\dim\widetilde{\mathcal{Y}}^{(\ell)}=i$ and $\dim\,\widetilde{\mathcal{U}}^{(\ell)} \ge i$
hold for each $\ell$. Moreover, consider an orthonormal matrix $\widetilde{U}\in\CC^{n \times i}$
with $\spanl\{\widetilde{U}\}\subseteq\widetilde{\mathcal{U}}^{(\ell)}$,
and an orthonormal basis matrix $\widetilde{Y}$ of $\widetilde{\mathcal{Y}}^{(\ell)}$,
and let $\mu(\cdot)$ be the Rayleigh quotient with respect to $H$.
If the matrix $\widetilde{R}=H\widetilde{Y}-\widetilde{Y}\widetilde{Y}^*H\widetilde{Y}$
has full rank, and
\[\|(H\widetilde{Y}-\widetilde{U}\widetilde{U}^*H\widetilde{Y})(\widetilde{R}^*\widetilde{R})^{-1/2}\|_2
 \le\widetilde{\gamma}<1\]
is fulfilled for each $\ell<L$ as in \eqref{prec2a}, then
\begin{equation}\label{thm2e1}
 \frac{\mu_i-\theta_i^{(L)}}{\theta_i^{(L)}-\mu_{s+1}}
 \le \left(\frac{\kappa_i+\widetilde{\gamma}(2-\kappa_i)}
 {(2-\kappa_i)+\widetilde{\gamma}\kappa_i}\right)^{2L}
 \frac{\mu_i-\theta_s^{(0)}}{\theta_s^{(0)}-\mu_{s+1}}
 \quad\mbox{with}\quad
 \kappa_i=\frac{\mu_{s+1}-\mu_n}{\mu_{i}-\mu_n}
\end{equation}
holds for the Ritz values produced by \eqref{bpg2}.
\item[(c)] If $\theta_s^{(0)}$ is located in $(\mu_{j+1},\mu_j]$ for a certain $j \ge s$,
consider the auxiliary iteration \eqref{bpg2a2}
using the same initial subspace $\spanl\{Y^{(0)}\}$ together with the invariant subspaces
$\widetilde{\mathcal{Z}}=\spanl\{z_{j-s+i+1},\ldots,z_j\}^{\perp}$
and $\widehat{\mathcal{Z}}=\spanl\{z_{j-s+1},\ldots,z_{j-s+i}\}^{\perp}$.
Assume for each $\ell$ that $\dim\spanl\{Y^{(\ell)}\}=s$,
$\dim\widetilde{\mathcal{Y}}^{(\ell)}=i$, and $\dim\widehat{\mathcal{Y}}^{(\ell)}=s{\,-\,}i$.
Then $\dim\,\widetilde{\mathcal{U}}^{(\ell)} \ge i$, and
a similar estimate for \eqref{bpg2} reads
\begin{equation}\label{thm2e2}
 \frac{\mu_{j-s+i}-\theta_i^{(L)}}{\theta_i^{(L)}-\mu_{j+1}}
 \le \left(\frac{\kappa_i+\widetilde{\gamma}(2-\kappa_i)}
 {(2-\kappa_i)+\widetilde{\gamma}\kappa_i}\right)^{2L}
 \frac{\mu_{j-s+i}-\theta_s^{(0)}}{\theta_s^{(0)}-\mu_{j+1}}
 \quad\mbox{with}\quad
 \kappa_i=\frac{\mu_{j+1}-\mu_n}{\mu_{j-s+i}-\mu_n}.
\end{equation}
\end{itemize}
\end{theorem}
\begin{proof}
(a) The trivial relation $\theta_i^{(\ell+1)}\ge\theta_i^{(\ell)}$
follows from the optimality of the Rayleigh-Ritz procedure.
For showing its strict version, we represent the trial subspace of \eqref{bpg2} by
\[\spanl\{Y^{(\ell)},U^{(\ell)}\}
 \quad\mbox{with}\quad
 U^{(\ell)}=Y^{(\ell)}\Theta_Y^{(\ell)}+N(HY^{(\ell)}-Y^{(\ell)}\Theta_Y^{(\ell)}).\]
If $\spanl\{Y^{(\ell)}\}$ contains no eigenvectors,
we use \cite[Lemma 3.1]{NEZ2019} where $U^{(\ell)}$
is analyzed within a BPG iteration with fixed step sizes. This implies
\[\theta_i^{(\ell+1)}=\theta_i(\spanl\{Y^{(\ell)},U^{(\ell)}\})
 \ge\theta_i(\spanl\{U^{(\ell)}\})>\theta_i(\spanl\{Y^{(\ell)}\})=\theta_i^{(\ell)}.\]

(b) Combining Lemma \ref{lmaux2} and Theorem \ref{thm1}
yields $\dim\widetilde{\mathcal{Y}}^{(\ell)}=i$, $\dim\,\widetilde{\mathcal{U}}^{(\ell)} \ge i$
and suggests a vector iteration concerning $\spanl\{\widetilde{y},\widetilde{u}\}$.
Moreover, for the respective smallest ($i$th largest) Ritz values
$\widetilde{\theta}_i^{(\ell)}$ and $\widetilde{\theta}_i^{(\ell+1)}$ of $H$
in $\widetilde{\mathcal{Y}}^{(\ell)}$ and $\widetilde{\mathcal{Y}}^{(\ell+1)}$,
the Courant-Fischer principles ensure
\[\mu_i\ge\widetilde{\theta}_i^{(\ell+1)}=\widetilde{\theta}'\ge\widetilde{\theta}^{\diamond}
 \ge\mu(\widetilde{y})\ge\widetilde{\theta}_i^{(\ell)}\]
(where $\widetilde{\theta}^{(\ell+1)}=\widetilde{\theta}'$ follows from Lemma \ref{lmaux2}),
i.e., the sequence $(\widetilde{\theta}_i^{(\ell)})_{\ell\in\mathbb{N}}$ is nondecreasing.
Thus $\widetilde{\theta}_i^{(\ell)}\ge\widetilde{\theta}_i^{(0)}\ge\theta_s^{(0)}>\mu_{s+1}$
holds so that the above Ritz values and $\mu(\widetilde{y})$
are all located in $(\mu_{s+1},\mu_i]$.
Then Lemma \ref{lmvit2} with $j{\,=\,}s$ leads to
\[\frac{\mu_i-\widetilde{\theta}^{\diamond}}{\widetilde{\theta}^{\diamond}-\mu_{s+1}}
 \le \left(\frac{\kappa_i+\widetilde{\gamma}(2-\kappa_i)}
 {(2-\kappa_i)+\widetilde{\gamma}\kappa_i}\right)^2
 \frac{\mu_i-\mu(\widetilde{y})}{\mu(\widetilde{y})-\mu_{s+1}}
 \quad\mbox{with}\quad
 \kappa_i=\frac{\mu_{s+1}-\mu_n}{\mu_{i}-\mu_n}\]
which can be extended as
\[\frac{\mu_i-\widetilde{\theta}_i^{(\ell+1)}}{\widetilde{\theta}_i^{(\ell+1)}-\mu_{s+1}}
 \le \left(\frac{\kappa_i+\widetilde{\gamma}(2-\kappa_i)}
 {(2-\kappa_i)+\widetilde{\gamma}\kappa_i}\right)^2
 \frac{\mu_i-\widetilde{\theta}_i^{(\ell)}}{\widetilde{\theta}_i^{(\ell)}-\mu_{s+1}}\]
by using the monotonicity of $(\mu_i-*)/(*-\mu_{s+1})$.
Recursively applying this intermediate estimate results in \eqref{thm2e1}
due to $\theta_i^{(L)}\ge\widetilde{\theta}_i^{(L)}$ and $\widetilde{\theta}_i^{(0)}\ge\theta_s^{(0)}$.

(c) Combining Lemma \ref{lmaux3} and Theorem \ref{thm1}
leads to $\dim\,\widetilde{\mathcal{U}}^{(\ell)} \ge i$ and a vector iteration.
The estimate \eqref{thm2e2} is trivial for $\theta_i^{(L)}\ge\mu_{j-s+i}$.
If $\theta_i^{(L)}<\mu_{j-s+i}$, we get
\[\mu_{j-s+i}>\theta_i^{(L)}\ge\widetilde{\theta}_i^{(L)}
 \ge\widetilde{\theta}_i^{(\ell+1)}=\widetilde{\theta}'\ge\widetilde{\theta}^{\diamond}
 \ge\mu(\widetilde{y})\ge\widetilde{\theta}_i^{(\ell)}
 \ge\widetilde{\theta}_i^{(0)}\ge\theta_s^{(0)}>\mu_{j+1}\]
for $\ell<L$ similarly to (b). Then \eqref{thm2e2} is derived by
Lemma \ref{lmvit2} with $j \ge s$ and a recursive reformulation
as well as monotonicity arguments.
\end{proof}

\begin{remark}\label{thmrem1}
Theorem \ref{thm2} extends the estimates for a BPG iteration with fixed step sizes
from \cite[Theorem 3.2 and Theorem 3.3]{NEZ2019} to the iteration \eqref{bpg2}
with implicitly optimized step sizes. The statement (a) indicates that the $i$th Ritz value
strictly increases until some eigenvectors are enclosed by the subspace iterate.
In addition, a reformulation of \cite[(3.7)]{NEZ2014} leads to the sharp estimate
\begin{equation}\label{thm2e3}
 \frac{\mu_j-\theta_s^{(\ell+1)}}{\theta_s^{(\ell+1)}-\mu_{j+1}}
 \le \left(\frac{\kappa+\gamma(2-\kappa)}
 {(2-\kappa)+\gamma\kappa}\right)^2
 \frac{\mu_j-\theta_s^{(\ell)}}{\theta_s^{(\ell)}-\mu_{j+1}}
 \quad\mbox{with}\quad
 \kappa=\frac{\mu_{j+1}-\mu_n}{\mu_j-\mu_n}
\end{equation}
for the $s$th Ritz value in the case $\theta_s^{(\ell)}\in(\mu_{j+1},\mu_j)$
with $j \ge s$ using the quality parameter $\gamma$ from \eqref{prec2}.
Combining this with (a) shows that $\theta_s^{(\ell)}$ can converge
to an eigenvalue $\mu_j$ with $j>s$ and otherwise can exceed $\mu_{s+1}$.
If $\theta_s^{(\ell)}>\mu_{s+1}$ occurs, we can reset the index $\ell$ as $0$
and apply the statement (b) to the further steps.
The statement (c) formally generalizes (b) to arbitrarily located $\theta_s^{(0)}$
and provides a supplement to \eqref{thm2e3} for discussing
the convergence of the $i$th Ritz value in the first steps of \eqref{bpg2}.
The assumption on subspace dimensions is usually fulfilled
in numerical tests with randomly generated initial guesses.
\end{remark}

\begin{remark}\label{thmrem2}
For evaluating the quality parameter $\widetilde{\gamma}$
in the estimates \eqref{thm2e1} and \eqref{thm2e2}, we can
follow the introduction of \eqref{prec2a} and determine
the auxiliary subspaces $\widetilde{\mathcal{Y}}^{(\ell)}=\spanl\{Y^{(\ell)}\}\cap\widetilde{\mathcal{Z}}$
and $\,\widetilde{\mathcal{U}}^{(\ell)}=\spanl\{U^{(\ell)}\}\cap\widetilde{\mathcal{Z}}$
via the invariant subspace $\spanl\{z_{j-s+1},\ldots,z_j\}$.
Furthermore, it is remarkable that \eqref{thm2e3} with $j=s$ implies
\[\frac{\mu_s-\theta_s^{(L)}}{\theta_s^{(L)}-\mu_{s+1}}
 \le \left(\frac{\kappa+\gamma(2-\kappa)}
 {(2-\kappa)+\gamma\kappa}\right)^{2L}
 \frac{\mu_s-\theta_s^{(0)}}{\theta_s^{(0)}-\mu_{s+1}}
 \quad\mbox{with}\quad
 \kappa=\frac{\mu_{s+1}-\mu_n}{\mu_s-\mu_n}\]
which is similar to \eqref{thm2e1} with $i=s$. These two estimates
for the $s$th Ritz value in the final phase of \eqref{bpg2}
only differ in the quality parameters $\gamma$ and $\widetilde{\gamma}$.
\end{remark}

\begin{remark}\label{thmrem3}
In comparison to the results from \cite{OVT2006},
our multi-step estimate \eqref{thm2e1} indicates that
the single-step convergence rate is asymptotically bounded by
$\widetilde{q}_i^2$ with $\widetilde{q}_i=\big(\kappa_i+\widetilde{\gamma}(2-\kappa_i)\big)/
\big((2-\kappa_i)+\widetilde{\gamma}\kappa_i\big)$ similarly to
the asymptotic convergence factor $q_{k,m}$ presented in \cite[Corollary 1]{OVT2006}
(despite typo with a redundant exponent $2$). If adapted to Theorem \ref{thm2}
(with $k \to i$ and $m \to s$), $q_{k,m}$ becomes
\[q_i=\frac{\kappa+\gamma(2-\kappa)}
 {(2-\kappa)+\gamma\kappa}
 \quad\mbox{with}\quad
 \kappa=\frac{\mu_{s+1}}{\mu_{i}}\]
which is slightly larger than $\widetilde{q}_i$ for $\widetilde{\gamma}=\gamma$.
However, the nonasymptotic estimate in \cite[Corollary 1]{OVT2006}
is formulated for a sum of Ritz value errors corresponding to
$\sum_{t=1}^i(\mu_t-\theta_t^{(\ell)})$. Therein the convergence bound
contains $q_i^2$ and $\sum_{t=1}^s(\mu_t-\theta_t^{(\ell)})$
together with a technical term which is not explicitly given.
The main estimate in \cite[Theorem 3]{OVT2006} uses a convergence factor
depending on certain angles and a ratio corresponding to
$\mu_{s+1}/\theta_i^{(\ell)}$ as a counterpart of the above
$\kappa=\mu_{s+1}/\mu_{i}$ (in the original formulation,
$\mu_k^i/\mu_{m+1}$ should be corrected as $\mu_{m+1}/\mu_k^i$).
In Theorem \ref{thm2}, we have achieved a concise convergence factor
by using the alternative quality parameter $\widetilde{\gamma}$.
The convergence rates of individual Ritz values do not need to
be analyzed in a mixed form.
\end{remark}

Finally, we reformulate Theorem \ref{thm2} as explicit statements
for the BPG iteration \eqref{bpg1} by using the substitutions
\eqref{subst} and \eqref{rpr}.

\begin{theorem}\label{thm3}
Consider the generalized eigenvalue problem \eqref{evp}
with $A$-orthonormal eigenvectors $w_1,\ldots,w_n$ of $(M,A)$
associated with the eigenvalues $\mu_1\ge\cdots\ge\mu_n$,
and let $\theta_1^{(\ell)}\ge\cdots\ge\theta_s^{(\ell)}$
be the Ritz values of $(M,A)$ in the subspace iterate $\spanl\{V^{(\ell)}\}$
of \eqref{bpg1}. Therein
$\Theta_V^{(\ell)}=\diag(\theta_1^{(\ell)},\ldots,\theta_s^{(\ell)})$,
$R_V^{(\ell)}=MV^{(\ell)}-AV^{(\ell)}\Theta_V^{(\ell)}$,
and the preconditioner $\widetilde{T}$ satisfies \eqref{prec1}.
Then the following statements hold:
\begin{itemize}
\setlength\itemsep{1ex}
\item[(a)] The Ritz values produced by \eqref{bpg1} fulfill $\theta_i^{(\ell+1)}\ge\theta_i^{(\ell)}$
for each $\ell$ and $i\in\{1,\ldots,s\}$. If there are no eigenvectors in $\spanl\{V^{(\ell)}\}$,
then $\theta_i^{(\ell+1)}>\theta_i^{(\ell)}$.
\item[(b)] If $\theta_s^{(0)}>\mu_{s+1}$, consider the auxiliary iteration
\begin{equation}\label{bpg2a3}
\begin{split}
 &\widetilde{\mathcal{V}}^{(\ell)}=\spanl\{V^{(\ell)}\}\cap\widetilde{\mathcal{W}},\qquad
 \widehat{\mathcal{V}}^{(\ell)}=\spanl\{V^{(\ell)}\}\cap\widehat{\mathcal{W}},\\[1ex]
 &U^{(\ell)}=V^{(\ell)}\Theta_V^{(\ell)}
 +T(MV^{(\ell)}-AV^{(\ell)}\Theta_V^{(\ell)})
 \quad\mbox{with}\quad T=\omega\widetilde{T} \ \ \mbox{from \eqref{prec1}},\\[1ex]
 &\,\widetilde{\mathcal{U}}^{(\ell)}=\spanl\{U^{(\ell)}\}\cap\widetilde{\mathcal{W}},\qquad
 \,\widehat{\mathcal{U}}^{(\ell)}=\spanl\{U^{(\ell)}\}\cap\widehat{\mathcal{W}},\\[1ex]
 &\widetilde{V}^{(\ell+\frac12)} \quad\xleftarrow{\mathrm{RR}[M,A,\widetilde{i}]}\quad
 \widetilde{\mathcal{V}}^{(\ell)}+\,\widetilde{\mathcal{U}}^{(\ell)},\qquad
 \widehat{V}^{(\ell+\frac12)} \quad\xleftarrow{\mathrm{RR}[M,A,\widehat{i}]}\quad
 \widehat{\mathcal{V}}^{(\ell)}+\,\widehat{\mathcal{U}}^{(\ell)},\\[1ex]
 &V^{(\ell+1)} \quad\xleftarrow{\mathrm{RR}[M,A,\widetilde{s}]}\quad
 \spanl\{\widetilde{V}^{(\ell+\frac12)},\widehat{V}^{(\ell+\frac12)}\}
\end{split}
\end{equation}
using the same initial subspace $\spanl\{V^{(0)}\}$ together with the invariant subspaces
$\widetilde{\mathcal{W}}=\spanl\{w_{i+1},\ldots,w_s\}^{\perp_A}$
and $\widehat{\mathcal{W}}=\spanl\{w_1,\ldots,w_i\}^{\perp_A}$.
Then $\widetilde{i}=\dim\widetilde{\mathcal{V}}^{(\ell)}=i$ and $\dim\,\widetilde{\mathcal{U}}^{(\ell)} \ge i$
hold for each $\ell$. Moreover, consider an $A$-orthonormal matrix $\widetilde{U}\in\CC^{n \times i}$
with $\spanl\{\widetilde{U}\}\subseteq\widetilde{\mathcal{U}}^{(\ell)}$,
and an $A$-orthonormal basis matrix $\widetilde{V}$ of $\widetilde{\mathcal{V}}^{(\ell)}$,
and let $\mu(\cdot)$ be the Rayleigh quotient with respect to $(M,A)$.
If $\widetilde{R}=A^{-1}M\widetilde{V}-\widetilde{V}\widetilde{V}^*M\widetilde{V}$
has full rank, and $\widehat{R}=(A^{-1}M\widetilde{V}-\widetilde{U}\widetilde{U}^*M\widetilde{V})
(\widetilde{R}^*A\widetilde{R})^{-1/2}$ fulfills
\,$\|\widehat{R}^*A\widehat{R}\|_2^{1/2}\le\widetilde{\gamma}<1$\,
for each $\ell<L$, then the multi-step estimate
\eqref{thm2e1} holds for the Ritz values produced by \eqref{bpg1}.
\item[(c)] If $\theta_s^{(0)}$ is located in $(\mu_{j+1},\mu_j]$ for a certain $j \ge s$,
consider the auxiliary iteration \eqref{bpg2a3}
using the same initial subspace $\spanl\{V^{(0)}\}$ together with the invariant subspaces
\,$\widetilde{\mathcal{W}}=\spanl\{w_{j-s+i+1},\ldots,w_j\}^{\perp_A}$\,
and \,$\widehat{\mathcal{W}}=\spanl\{w_{j-s+1},\ldots,w_{j-s+i}\}^{\perp_A}$.
Assume for each $\ell$ that $\dim\spanl\{V^{(\ell)}\}=s$,
$\widetilde{i}=\dim\widetilde{\mathcal{V}}^{(\ell)}=i$,
and $\widehat{i}=\dim\widehat{\mathcal{V}}^{(\ell)}=s{\,-\,}i$.
Then $\dim\,\widetilde{\mathcal{U}}^{(\ell)} \ge i$, and
a similar estimate for \eqref{bpg1} is given by \eqref{thm2e2}.
\end{itemize}
\end{theorem}

A further reformulation of Theorem \ref{thm3} for extending the convergence analysis
of the block preconditioned steepest descent iteration from \cite{NEZ2014} refers to
the computation of the smallest eigenvalues of the matrix pair $(A,M)$
for Hermitian positive definite $A,M\in\CC^{n \times n}$.
Therein the estimate \eqref{thm2e2} turns into
\begin{equation}\label{bpsde}
 \frac{\vartheta_i^{(L)}-\la_{j-s+i}}{\la_{j+1}-\vartheta_i^{(L)}}
 \le \left(\frac{\kappa_i+\widetilde{\gamma}(2-\kappa_i)}
 {(2-\kappa_i)+\widetilde{\gamma}\kappa_i}\right)^{2L}
 \frac{\vartheta_s^{(0)}-\la_{j-s+i}}{\la_{j+1}-\vartheta_s^{(0)}}
 \quad\mbox{with}\quad
 \kappa_i=\frac{\la_{j-s+i}(\la_n-\la_{j+1})}{\la_{j+1}(\la_n-\la_{j-s+i})}
\end{equation}
where $\la$ and $\vartheta$ denote eigenvalues and Ritz values
of $(A,M)$ in ascending order.

\begin{remark}\label{bpsd}
Although the parameter $\widetilde{\gamma}$ cannot easily be replaced by
$\gamma$ from \eqref{prec1} in our analysis,
mainly due to additional modifications of the preconditioned term $U^{(\ell)}$
by intersections in \eqref{bpg2a3}, the corresponding estimates
with $\gamma$ still provide reasonable bounds in numerical experiments.
An analysis directly using $\gamma$ should avoid additional modifications
as in the following auxiliary iteration:
\begin{equation}\label{bpg2b}
\begin{split}
 &\spanl\{\widetilde{V}^{(\ell)}\}=\spanl\{V^{(\ell)}\}\cap\widetilde{\mathcal{W}},\qquad
 \widetilde{V}^{(\ell+\frac12)} \quad\xleftarrow{\mathrm{RR}[M,A,\widetilde{i}]}\quad
 \spanl\{\widetilde{V}^{(\ell)},TR_{\widetilde{V}}^{(\ell)}\},\\[1ex]
 & \spanl\{\widehat{V}^{(\ell)}\}=\spanl\{V^{(\ell)}\}\cap\widehat{\mathcal{W}},\qquad
 \widehat{V}^{(\ell+\frac12)} \quad\xleftarrow{\mathrm{RR}[M,A,\widehat{i}]}\quad
 \spanl\{\widehat{V}^{(\ell)},TR_{\widehat{V}}^{(\ell)}\},\\[1ex]
 &V^{(\ell+1)} \quad\xleftarrow{\mathrm{RR}[M,A,\widetilde{s}]}\quad
 \spanl\{\widetilde{V}^{(\ell+\frac12)},\widehat{V}^{(\ell+\frac12)}\}
\end{split}
\end{equation}
where $R$ denotes block residuals.
In the case $\widetilde{i}=\dim\spanl\{\widetilde{V}^{(\ell)}\}=1$, the first partial iteration
in \eqref{bpg2b} is a vectorial gradient iteration. A geometric relation
between two successive iterates $\widetilde{v}$ and $\widetilde{v}'$
can be derived based on \cite[Theorem 3.1]{NEY2012}, namely,
there is a rational function $f(\cdot)$ satisfying
$\widetilde{v}^{\diamond}=f(A^{-1}M)\widetilde{v}$ and $\mu(\widetilde{v}')\ge\mu(\widetilde{v}^{\diamond})$.
Moreover, $\widetilde{v}^{\diamond}$ can be regarded as the next iterate
generated by a special preconditioner $T^{\diamond}$.
Therefore the convergence of the largest Ritz value in \eqref{bpg2b}
is decelerated by using such $T^{\diamond}$.
The corresponding first partial iteration can be simplified 
since $\spanl\{\widetilde{V}^{(\ell+\frac12)}\}\subseteq\widetilde{\mathcal{W}}$
is ensured by $\widetilde{v}^{\diamond}=f(A^{-1}M)\widetilde{v}$.
This results in \eqref{thm2e2} with $\widetilde{\gamma}=\gamma$ for $i=1$.
Nevertheless, a generalization to arbitrary $i\in\{1,\ldots,s\}$
requires further assumptions on partial iterations.
Occasionally, we can apply the estimate with $\gamma$ for $i=1$ similarly to a deflation,
i.e., analyzing the convergence rate of the $(i{\,+\,}1)$th Ritz value
provided that the first $i$ Ritz values are sufficiently close to the target eigenvalues.
\end{remark}

\section{Numerical experiments}

We consider several numerical examples in order to demonstrate the main results
and discuss their accuracy. In the first example, we implement the accompanying
iteration \eqref{bpg2} for a test matrix from \cite{NEZ2019},
and illustrate Theorem \ref{thm2}. The further examples using
discretized Laplacian eigenvalue problems are concerned with
the BPG iteration \eqref{bpg1} and Theorem \ref{thm3}.

\textbf{Example I.}
We reuse the diagonal matrix $H=\diag(\mu_1,\ldots,\mu_n)$
from \cite[Experiment I]{NEZ2019} with $n=6000$
and \,$\mu_i=10.07-0.01\,i$\, for $i\le 6$.
The further eigenvalues (diagonal entries) of $H$
are given by equidistant points between $9$ and $1$.

\begin{figure}[htbp]
\begin{center}
\includegraphics[width=\textwidth]{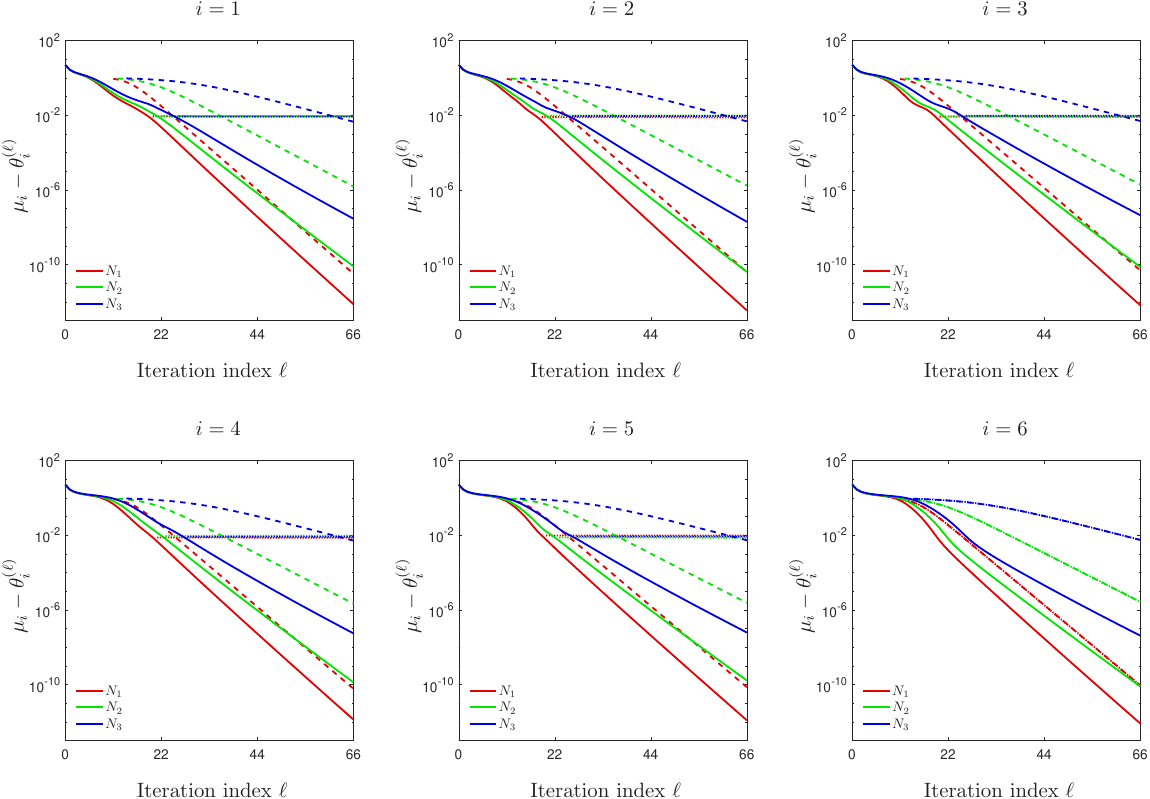}\quad
\end{center}
\par\vskip -2ex
\caption{\small Cluster robustness of the accompanying iteration \eqref{bpg2}
 applied to Example I. Solid curves: Ritz value errors in the slowest run
 among $1000$ runs with random initial subspaces.
 Dashed curves: Bounds determined by Theorem \ref{thm2}.
 Dotted curves: Bounds based on single-step estimates with neighboring eigenvalues.}
\label{fig1}
\end{figure}

We implement the iteration \eqref{bpg2} with the block size $s=6$
where the target eigenvalues $\mu_1,\ldots,\mu_s$ are tightly clustered.
We test three preconditioners, denoted by $N_1,\,N_2,\,N_3$.
The first one is simply $N_1=I$, whereas $N_2$ and $N_3$ are generated
by random sparse perturbations of $I$, namely,
\,\verb|N=eta*sprand(n,n,5/n); N=N'+I+N|\,
with $\eta\in\{0.09,\,0.16\}$. For each preconditioner,
we compare $1000$ runs with random initial subspaces,
and illustrate the slowest run with respect to
the Ritz value errors $\mu_i-\theta_i^{(\ell)}$, $i\in\{1,\ldots,6\}$
by solid curves in Figure \ref{fig1}.
This immediately reflects the monotone convergence in the statement (a)
of Theorem \ref{thm2}. 

For checking the statement (b),
we determine the quality parameter $\widetilde{\gamma}$
by evaluating \eqref{prec2a} within the auxiliary iteration \eqref{bpg2a2}
for each iteration step after the Ritz value $\theta_s^{(\ell)}$ exceeds $\mu_{s+1}$.
The corresponding maximum is used as $\widetilde{\gamma}$
in the estimate \eqref{thm2e1} with an index adaptation. Therein
\[\widetilde{\gamma}=0 \ \mbox{for} \ N_1, \quad
 \widetilde{\gamma}\approx 0.2429 \ \mbox{for} \ N_2, \quad
 \widetilde{\gamma}\approx 0.5285 \ \mbox{for} \ N_3.\]
The resulting bounds for $\mu_i-\theta_i^{(\ell)}$
are plotted by dashed curves in Figure \ref{fig1}. 
In addition, their counterparts based on the single-step estimates
from \cite{NEZ2014}, i.e., those using
$(\mu_{i+1}-\mu_n)/(\mu_{i}-\mu_n)$
instead of $(\mu_{s+1}-\mu_n)/(\mu_{i}-\mu_n)$ in \eqref{thm2e1},
are displayed by dotted curves.

These two types of bounds coincide for $i=s=6$.
The difference between them is substantial for $i\in\{1,\ldots,5\}$
due to $\mu_{i}\approx\mu_{i+1}$. The bounds in dashed curves
clearly reflect the cluster robustness, whereas the bounds
in dotted curves wrongly predict a stagnation.
Furthermore, the accuracy of bounds in dashed curves
apparently depends on the accuracy of preconditioning,
and could be improved for less accurate preconditioners. This motivates
a future task for defining a more effective quality parameter.

The statement (c) can be checked in a similar way. We omit
the illustration since it only concerns a few iteration steps
for random initial subspaces. A reasonable illustration requires
certain special initial subspaces.

\textbf{Example II.}
We consider the Laplacian eigenvalue problem on the rectangle domain
$[0,\,2]\times[0,\,1]$ with a slit $\{1\}\times[0.1,\,0.9]$
and homogeneous Dirichlet boundary conditions.
The five-point star discretization with the mesh size
$1/70$ results in a standard eigenvalue problem
which can be reformulated as \eqref{evp} for $n=9534$ and $M=I$.
The six largest eigenvalues build two tight clusters
$\{\mu_1,\mu_2\}$ and $\{\mu_3,\ldots,\mu_6\}$.

\begin{figure}[htbp]
\begin{center}
\includegraphics[width=\textwidth]{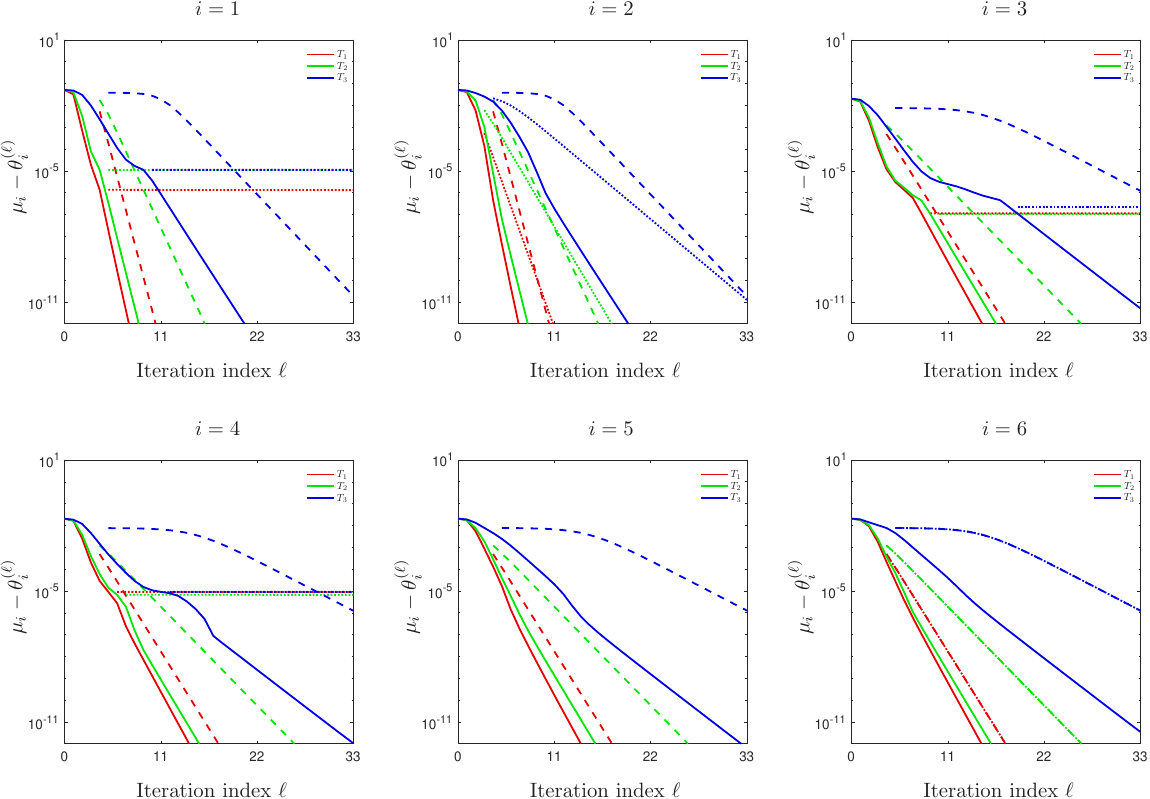}\quad
\end{center}
\par\vskip -2ex
\caption{\small Cluster robustness of the BPG iteration \eqref{bpg1}
 applied to Example II. Solid curves: Ritz value errors in the slowest run
 among $1000$ runs with random initial subspaces.
 Dashed curves: Bounds determined by Theorem \ref{thm3}.
 Dotted curves: Bounds based on single-step estimates with neighboring eigenvalues.}
\label{fig2}
\end{figure}

The BPG iteration \eqref{bpg1} with the block size $s=6$
is implemented for three preconditioners $T_1,\,T_2,\,T_3$ constructed by
\begin{center}
\,\verb|ichol(A,struct('type','ict','droptol',eta))|\,
for $\eta\in\{10^{-5},\,10^{-4},\,10^{-3}\}$.
\end{center}
Similarly to Example I, Ritz value errors in the slowest run
concerning $1000$ random initial subspaces
are illustrated by solid curves in Figure \ref{fig2}.

We particularly demonstrate the statement (b) of Theorem \ref{thm3}.
Therein the quality parameter $\widetilde{\gamma}$ is determined
for each iteration step after $\theta_s^{(\ell)}>\mu_{s+1}$
by using the auxiliary iteration \eqref{bpg2a3}.
The respective maxima are
\[\widetilde{\gamma}\approx 0.0295  \ \mbox{for} \ T_1, \quad
 \widetilde{\gamma}\approx 0.2344 \ \mbox{for} \ T_2, \quad
 \widetilde{\gamma}\approx 0.5168  \ \mbox{for} \ T_3.\]
The resulting bounds in dashed curves are appropriate
for each Ritz value. Their counterparts based on \cite{NEZ2014}
in dotted curves are only reasonable for $i\in\{2,\,6\}$
where $\mu_{i}$ and $\mu_{i+1}$ are not clustered
and the bounds are more accurate in the first steps.
Moreover, the dotted curves cannot be drawn for $i=5$
since $\mu_5$ and $\mu_6$ coincide.
The new bounds are thus advantageous
for clustered and multiple eigenvalues.

\begin{figure}[htbp]
\mbox{\quad}\vspace{6mm}
\begin{center}
\includegraphics[width=\textwidth]{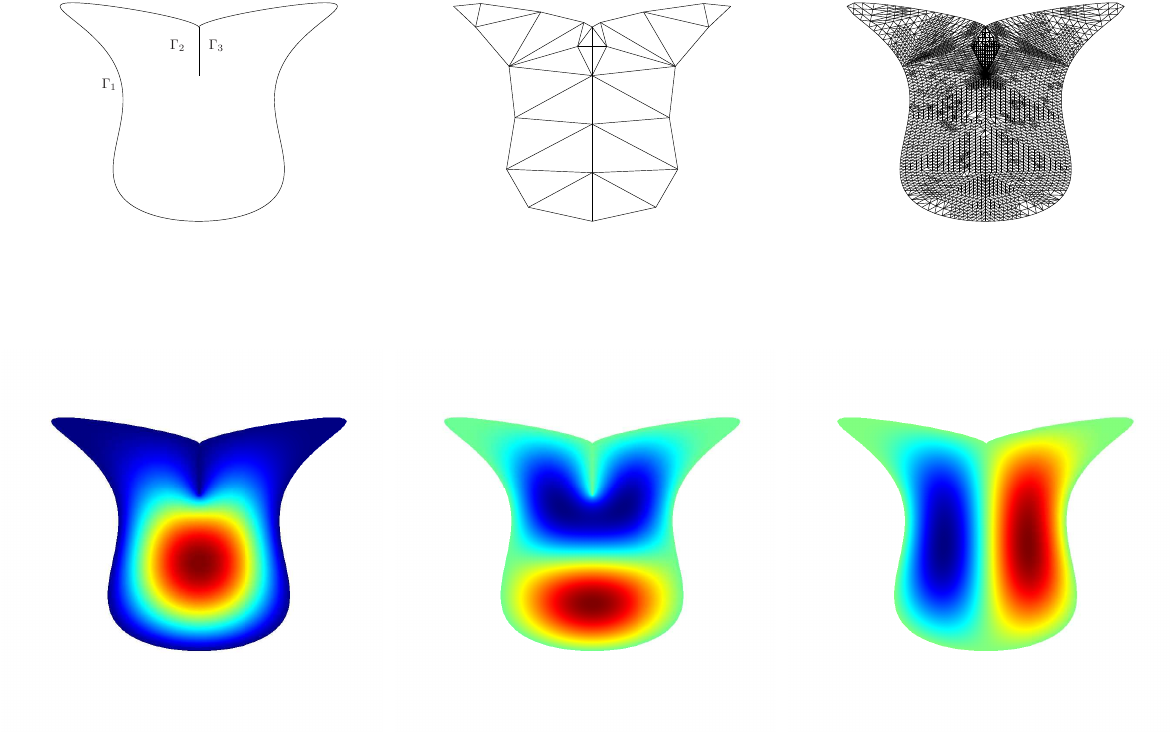}\quad
\end{center}
\par\vskip -2ex
\caption{\small Laplacian eigenvalue problem for Example III.
First row: domain, initial grid and an adaptively refined grid.
Second row: approximate eigenfunctions (top view)
associated with the three smallest operator eigenvalues
whose residuals are used for the grid refinement.}
\label{figdm}
\end{figure}

\textbf{Example III.}
We consider the Laplacian eigenvalue problem on a 2D tulip-like domain
with homogeneous Dirichlet boundary conditions; see Figure \ref{figdm}.
The boundary consists of three parts:
\[\begin{split}
 &\Gamma_1=\big\{\big(1.2\sin(t){\,+\,}0.3\sin(4t),\
  {\,-\,}\cos(t){\,-\,}0.5\cos(2t)\big)^T;\ t\in[-\pi,\,\pi)\big\},\\
 &\Gamma_2=\big\{\big(0,\ 0.5\,t\big)^T;\ t\in(0,\,1)\big\},\qquad
  \Gamma_3=\big\{\big(0,\ 0.5(1-t)\big)^T;\ t\in(0,\,1]\big\}.
\end{split}\]
We generate matrix eigenvalue problems successively
by an adaptive finite element discretization depending on
residuals of approximate eigenfunctions associated with
the three smallest operator eigenvalues; cf.~\cite[Appendix]{NEZ2019}
and some relevant graphics in Figure \ref{figdm}.
We repeat the numerical experiments in Example II
for the matrix pair $(M,A)$ from the $41$st grid of the discretization
with $n=1{,}522{,}640$ degrees of freedom.
The largest eigenvalues of $(M,A)$ approximate the reciprocals
of the smallest operator eigenvalues.

We observe again the BPG iteration \eqref{bpg1} with the block size $s=6$.
The target eigenvalues are partially clustered ($\mu_2\approx\mu_3$).
Concerning the comparison between new results in Theorem \ref{thm3}
and previous results based on \cite{NEZ2014}, we first compare
their decisive terms $\kappa_i=(\mu_{s+1}-\mu_n)/(\mu_{i}-\mu_n)$
and $\widehat{\kappa}_i=(\mu_{i+1}-\mu_n)/(\mu_{i}-\mu_n)$.
For $i\in\{1,\ldots,6\}$, we have
\begin{equation}\label{kc}
\begin{split}
 &\kappa_i\in\{0.1993,\,0.4628,\,0.4687,\,0.6279,\,0.7307,\,0.9459\},\\
 &\widehat{\kappa}_i\in\{0.4306,\,0.9875,\,0.7464,\,0.8593,\,0.7725,\,0.9459\}.
\end{split}
\end{equation}
The test preconditioners $T_1,\,T_2,\,T_3$ are constructed by
\,\verb|ichol|\, with $\eta\in\{10^{-7},\,10^{-6},\,10^{-5}\}$
as \,\verb|droptol|. The quality parameter $\widetilde{\gamma}$
with respect to the auxiliary iteration \eqref{bpg2a3} reads
\[\widetilde{\gamma}\approx 0.1116  \ \mbox{for} \ T_1, \quad
 \widetilde{\gamma}\approx 0.3563 \ \mbox{for} \ T_2, \quad
 \widetilde{\gamma}\approx 0.6652  \ \mbox{for} \ T_3.\]

\begin{figure}[htbp]
\begin{center}
\includegraphics[width=\textwidth]{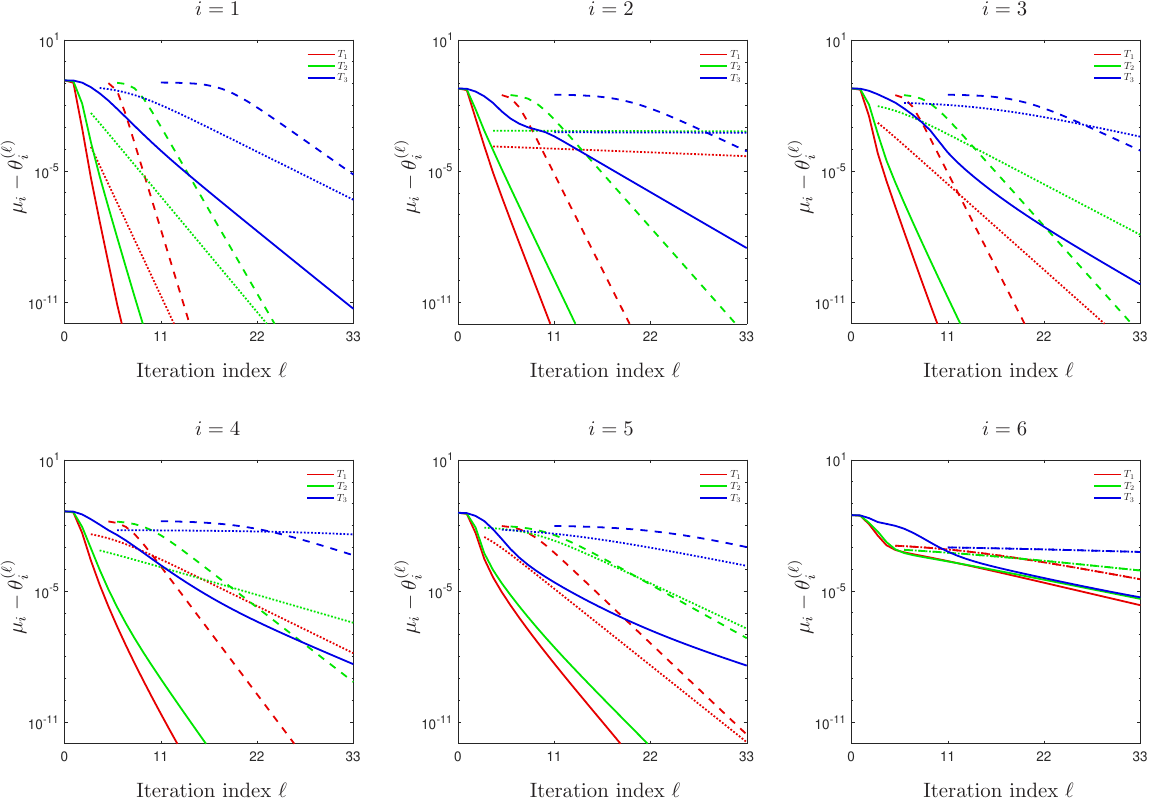}\quad
\end{center}
\par\vskip -2ex
\caption{\small Cluster robustness of the BPG iteration \eqref{bpg1}
 applied to Example III. Solid curves: Ritz value errors in the slowest run
 among $1000$ runs with random initial subspaces.
 Dashed curves: Bounds determined by Theorem \ref{thm3}.
 Dotted curves: Bounds based on single-step estimates with neighboring eigenvalues.}
\label{fig3}
\end{figure}

Figure \ref{fig3} presents a bound comparison for Ritz value errors
in the slowest run concerning $1000$ random initial subspaces.
The dashed curves display new bounds from Theorem \ref{thm3}.
They generally have steeper slopes than
the dotted curves containing bounds based on \cite{NEZ2014}.
The slope difference mainly depends on the terms
$\kappa_i$ and $\widehat{\kappa}_i$; cf.~their values given in \eqref{kc}.
The maximal difference appears for $i=2$ where
the dotted curves are almost constant. As an explanation, we note that
the corresponding $\kappa$ value
$\widehat{\kappa}_2\approx0.9875$ is close to $1$, and the convergence factor
is at least $\widehat{\kappa}_2/(2-\widehat{\kappa}_2)$
for each test preconditioner. Such an overestimation can also be caused by
a slightly smaller $\kappa$ value between $0.85$ and $0.95$ for a moderate
preconditioner; cf.~the blue curves for $i\in\{4,\,6\}$ corresponding to $T_3$
combined with $\widehat{\kappa}_4$ (dotted),
$\kappa_6$ and $\widehat{\kappa}_6$ (dashed and dotted).
Deriving sharper bounds in the case of moderate preconditioners
is potentially important for large-scale discretized eigenvalue problems
where generating more accurate preconditioners,
e.g. with $\widetilde{\gamma}<0.5$, is costly with respect to inner steps
and the total time.

\section*{Conclusion}
The cluster robustness of block preconditioned gradient (BPG) eigensolvers
with sufficiently large block sizes is studied by deriving proper convergence bounds
of individual Ritz values. A basic argument
in our analysis is that the Rayleigh-Ritz (RR) approximation
in the trial subspace of BPG can be decelerated by applying RR to certain
lower-dimensional subspaces. This motivates auxiliary iterations
whose iterates are orthogonal to eigenvectors associated with some
possibly clustered eigenvalues. The relevant eigenvalues in the resulting bound
are thus not close to each other and reflect a cluster-independent convergence rate.
The construction of such auxiliary iterations is relatively easy for exact-inverse preconditioning
by using the classical analysis of an abstract block iteration \cite{KNY1987}.
The previous analysis \cite{NEZ2019} deals with
an arbitrary Hermitian positive definite preconditioner, but focuses on fixed step sizes
which correspond to the block power method rather than a block gradient iteration.
Therein an alternative quality parameter for the preconditioner leads to
concise bounds under weaker assumptions in comparison to \cite{BPK1996,OVT2006}.
This approach is upgraded in the present paper by adapting some geometric arguments
from our analysis of the (block) preconditioned steepest descent iteration
\cite{NEY2012,NEZ2014}. The achieved multi-step estimates improve
the sumwise estimates from \cite{OVT2006} in the sense of
more intuitive convergence factors and the applicability to individual Ritz values.
It is remarkable that BPG as two-block iterations are not necessarily cluster robust
for small block sizes. This drawback can be overcome by three(or more)-block iterations
such as LOBPCG and restarted Davidson methods \cite{STA2007,SMC2007,WXS2019}.
Extending our analysis of BPG to more powerful eigensolvers
is desirable in our future research.

\end{document}